\documentclass[12pt]{article}
\usepackage{latexsym}
\usepackage{amsmath,a4wide, amssymb}
\usepackage[dutch]{babel}

\usepackage{mathrsfs}
\usepackage[T1]{fontenc}
\usepackage{verbatim} 
\usepackage{holtpolt}
\usepackage{fancyheadings} 

\newcommand{\dk}[1]{\mathbf{#1}}
\newcommand{\gt}[1]{\gothic{#1}}
\newcommand{\dg}[1]{\boldsymbol{#1}}

\newcommand{\mr}[1]{\mathrm{#1}}
\newcommand{\ms}[1]{\mathscr{#1}}

\newcommand{\itl}[1]{\mathit{#1}}
\newcommand{\msf}[1]{\mathsf{#1}}
\newcommand{\B}[1]{\Bbb{#1}}
\DeclareMathAlphabet\gothic{U}{euf}{m}{n}

\newcommand{\dR}{\B{R}\,}
\newcommand{\dC}{\B{C}}
\newcommand{\dL}{\B{L}}

\newcommand{\rd}{{\,\rm d}}

\makeatletter
\@addtoreset{equation}{section}
\makeatother

\DeclareMathAlphabet\gothic{U}{euf}{m}{n}

\newcommand{\osx}{\underline{x}}
\newcommand{\osy}{\underline{y}}
\newcommand{\osa}{\underline{a}}
\newcommand{\osb}{\underline{b}}
\newcommand{\ose}{\underline{e}}
\newcommand{\osz}{\underline{z}}
\newcommand{\osw}{\underline{w}}

\newcommand{\cm}[1]{\mathcal{#1}}
\newcommand{\os}[1]{\underline{#1}}

\newcommand{\pa}{\partial}
\newcommand{\col}{\mathrm{col}}
\newcommand{\sC}{\mathscr{C}}

\newcommand{\sF}{\mathscr{F}}
\newcommand{\ssF}{\widetilde{\sF}}

\newcommand{\sQ}{\mathscr{Q}}

\newcommand{\ssG}{\widetilde{\sG}} 
\newcommand{\sL}{\mathscr{L}}

\newcommand{\sG}{\mathscr{G}}
\newcommand{\sO}{\mathscr{O}}
\newcommand{\sV}{\mathscr{V}}

\newcommand{\rK}{\mr{K}}

\newcommand{\rL}{\mr{L}}
\newcommand{\rI}{\mr{I}}

\newtheorem{lemma}{Lemma}[section] 
\newtheorem{thm}[lemma]{Theorem}  
\newtheorem{cor}[lemma]{Corollary}
\newtheorem{voorb}[lemma]{Example} 
\newtheorem{voorbn}[lemma]{Examples }
\newtheorem{rem}[lemma]{Remark}
\newtheorem{definitie}[lemma]{Definition}
\newtheorem{rems}[lemma]{Remarks}

\newtheorem{eig}[lemma]{Properties}
\newtheorem{cond}[lemma]{Condition}
\newtheorem{nnotation}[lemma]{Notation}
\newcommand{\eindebewijs}{\hfill{$\blacksquare$}}

\newtheorem{Nnotation}[lemma]{Notation and Definitions}

\newtheorem{nnotations}[lemma]{Notation}
\newtheorem{nnpotations}[lemma]{Notation and Properties}

\newenvironment{remarkn}{\begin{rem} \rm}{\end{rem}} 
\newenvironment{remarks}{\begin{rems} \rm \ \\}{\end{rems}} 
\newenvironment{examp}{\begin{voorb} \rm}{\end{voorb}}
\newenvironment{examps}{\begin{voorbn} \rm  }{\end{voorbn}}
\newenvironment{defn}{\begin{definitie} \rm}{\end{definitie}}
\newenvironment{notation}{\begin{nnotation} \rm}{\end{nnotation}}

\newenvironment{voorw}{\begin{cond} \rm}{\end{cond}} 

\newcommand{\be}{\begin{equation}}
\newcommand{\ee}{\end{equation}}
\newcommand{\bt}{\begin{thm}{\mbox{}\\}}
\newcommand{\et}{\end{thm}}
\newcommand{\bc}{\begin{voorw}{\mbox{}\\}}
\newcommand{\ec}{\end{voorw}}

\newcommand{\mA}{\mr{A}}
\newcommand{\mB}{\mr{B}}
\newcommand{\mG}{\mr{G}}
\newcommand{\mH}{\mr{H}}
\newcommand{\mK}{\mr{K}}
\newcommand{\mM}{\mr{M}}
\newcommand{\mN}{\mr{N}}
\newcommand{\mO}{\mr{O}}
\newcommand{\mP}{\mr{P}}
\newcommand{\mQ}{\mr{Q}}
\newcommand{\mR}{\mr{R}}
\newcommand{\mS}{\mr{S}}

\newcommand{\mX}{\mr{X}}
\newcommand{\mY}{\mr{Y}}
\newcommand{\mZ}{\mr{Z}}
\newcommand{\mW}{\mr{W}}

\newcommand{\dpsi}{\boldsymbol\Psi}
\newcommand{\hpsi}{\hat{\boldsymbol\Psi}}
\newcommand{\dA}{\mathbf{A}}

\newcommand{\dZ}{\mathbf{Z}}

\newcommand{\cA}{{\cm{A}}}

\newcommand{\Der}{{\cm{D}}} 
\newcommand{\cF}{{\cm{F}}}
\newcommand{\cH}{{\cm{H}}}

\newcommand{\cU}{\cm{U}}
\newcommand{\cV}{\cm{V}}

\newcommand{\cG}{\cm{G}}

\newcommand{\cL}{\cm{L}}
\newcommand{\hA}{\hat{\cm{A}}}

\newcommand{\hF}{\hat{\cm{F}}}
\newcommand{\hf}{\hat{f}}
\newcommand{\Tr}{\mathrm{Tr}} 
\newcommand{\Trc}[1]{\Tr\big\{#1\big\}}

\renewcommand{\Re}{\mathrm{Re}\,}
\renewcommand{\Im}{\mathrm{Im}\,}
\newcommand{\da}{\dagger}

\newcommand{\dint}{\displaystyle{\int}}
\newcommand{\dpa}{\dg{\pa}}

\newcommand{\sgn}{\mathrm{sgn}}
\newcommand{\ri}{\,\mathrm{i}}
\newcommand{\onul}{\underline{0}}
\newcommand{\bs}[1]{\overline{#1}}
\newcommand{\dhr}{\!\vartriangleleft\!} 

\newcommand{\gtG}{\gt{G}}
\newcommand{\gtg}{{\raisebox{0.5mm}{$\gt{g}$}}}
\newcommand{\fG}{\ms{C}^\infty(\dR^N\!:\!\gt{G})}
\newcommand{\fg}{\ms{C}^\infty(\dR^N\!:\!\gtg)}   

\newcommand{\fpsi}{\ms{C}^\infty(\dR^N\!:\!\dC^{r\ti c})}
\newcommand{\paf}[2]{\dfrac{\pa{\,#1}}{\pa{#2}}}           
\newcommand{\paft}[2]{\dfrac{\pa^2{#1}}{\pa{#2}^2}}
\newcommand{\gaf}[2]{\dfrac{\rd{#1}}{\rd{#2}}} 
\renewcommand{\div}{\mathrm{div}}
\newcommand{\grad}{\mathrm{grad}}
\newcommand{\rot}{\mathrm{rot}}
\renewcommand{\col}{\mathrm{col}}
\newcommand{\row}{\mathrm{row}}
\newcommand{\lra}{\leftrightarrow}
\newcommand{\proj}{\ms{P}_{\!\gt{g}}}

\newcommand{\al}{\alpha}
\newcommand{\bet}{\beta}
\newcommand{\ga}{\gamma}
\newcommand{\Ga}{\Gamma}
\newcommand{\del}{\delta}

\newcommand{\vep}{\varepsilon}

\renewcommand{\th}{\theta} 
\newcommand{\Th}{\Theta}

\newcommand{\ka}{\kappa}
\newcommand{\la}{\lambda}

\newcommand{\ti}{\times}
\newcommand{\glh}{\Big[}  
\newcommand{\grh}{\Big]}
\newcommand{\llh}{\dg{\{}} 
\newcommand{\rlh}{\dg{\}}}
\newcommand{\lk}{\dg{\,,\,}} 

\newcommand{\ths}{\theta^\star}
\newcommand{\rhos}{\rho^\star}

\newcommand{\gU}{\mathsf{U}}

\newcommand{\aA}{\mathcal{A}}
\newcommand{\aB}{\mathcal{B}}
\newcommand{\aF}{\mathcal{F}}
\newcommand{\aH}{\mathcal{H}}

\newcommand{\ad}{\mathrm{ad}}

\newcommand{\oa}{\os{a}}

\newcommand{\ox}{\os{x}}
\newcommand{\oA}{\os{A}}
\newcommand{\opsi}{\os{\psi}}

\newcommand{\vwL}{\msf{L}}
\newcommand{\vwG}{\msf{G}}

\begin{document} 
\  \vspace{-3cm}
\begin{center}
{\Huge\bf \
\\[3mm] Matrix Gauge Fields \\[6mm] and \\[8mm] Noether's Theorem}\\[12mm]
\end{center}
\begin{center}
{\large\bf J. de GRAAF}\\[6mm]
{ \bf Eindhoven University of Technology, Mathematics,}\\[3mm] 
{\bf Casa Reports 14-14, May 2014}
\
\\[10mm]
\end{center}

\vspace{12mm}
\begin{center}
{\bf Preface and Summary}
\end{center}

These notes are about  systems of 1st and 2nd order (non-)linear partial differential equations which are formed from a {\em Lagrangian density} $\vwL_\psi\,:\,\dR^N\to\dC~$,
\[ \text{Symbolically}:~~~~ \ox\mapsto\vwL_\psi(\ox)=\vwL(\opsi(\ox)\,;\nabla\opsi(\ox)\,;\ox)~,
\]
 by means of the usual Euler-Lagrange variational rituals.
 The non subscripted  $\vwL$ will denote the '{\em proto-Lagrangian}', which is a function of a finite number of  variables:
 \[
 \vwL~:~\dC^{r\ti c}\ti\dC^{Nr\ti c}\ti\dR^N~\to~\dC\,.
 \]
 In this $\vwL$  one has  to substitute matrix-valued functions  $\opsi:\dR^N\to\dC^{r\ti c}$ and
 $\nabla\opsi:\dR^N\to\dC^{Nr\ti c}$  for obtaining the Lagrangian density
$\vwL_\psi$. In our considerations the role and the special properties of the proto-Lagrangian $\vwL$ are crucial.\\[1mm]
These notes have been triggered by physicist's considerations: {\bf (1)} on obtaining the 'classical', that is the 'pre-quantized', wave equations for matter fields from variational principles, {\bf (2)} on conservation laws and {\bf (3)} on 'gauge field extensions'. For the humble mathematical anthropologist the rituals in physics textbooks have not much changed  during the last four decades.
Neither have they become much clearer. Compare e.g.
 [DM] and [W].\\[1mm]
The underlying notes give special attention  to the following \vspace{-2mm}
\begin{itemize}
\item In expressions (='equations') for Lagrange densities often  both $\opsi$ and its hermitean transposed $\opsi^\da$ appear. Are they meant as independent variables or not? Mostly, from the context the suggestion arises that 'variation' of $\opsi$ and 'variation' of $\opsi^\da$ lead to the same Euler-Lagrange equations. Why? Our remedy is  doubling the matrix entries in the proto-Lagrangian and thereby making  the Lagrangian density explicitly dependent on both $\opsi\,,\opsi^\da$ and their derivatives: So for $\vwL_\psi(\ox)$ we take expressions like $\sL_\psi(\ox)=\sL(\opsi(\ox)\,;\opsi(\ox)^\da\,;\nabla\opsi(\ox)\,;\nabla\opsi(\ox)^\da\,;\ox)$.  A suitable condition is then that the {\em Lagrangian functional}
\[\cL[\opsi]=\int_{\dR^N} \sL_\psi(\ox)\,\rd\ox
\]
only takes real values (Thm 2.4).
\item For 'free gauge fields' the situation is somewhat different. Now the dependent variables, named $\aA_\mu\,,1\leq\mu\leq N$, take their values in some fixed Lie-algebra $\gtg\subset\dC^{c\ti c}$. Although $\gtg$ mostly contains
complex matrices it is  a {\em real vector space} in interesting cases. (Note that $\gothic{u}(1)=
\ri\dR$ is a {\em real} vector space!). Therefore it needs a separate treatment.
\item The traditional conservation laws for quantities like energy, momentum, moment of momentum, $\ldots\,$, turn out to be  based  on {\em External Infinitesimal Symmetries}  of the proto-Lagrangian.
This means the existence of a couple of linear mappings

$\rK:\dC^{r\ti c}\to\dC^{r\ti c}\,,\rL:\dC^{Nr\ti c}\to\dC^{Nr\ti c}\,$, together with an affine mapping

$\ox\mapsto -s\oa+e^{sA}\ox~$, such that for all matrices $\mP\in\dC^{r\ti c}\,,\os{\mQ}\in\dC^{Nr\ti c}~$ and$~\ox\in\dR^N~$,
\[
\vwL(e^{s\rK}\mP\,;e^{s\rL}\os{\mQ}\,;-s\oa+e^{sA}\ox)\,=\,\vwL(\mP\,,\os{\mQ}\,;\ox)\,+\,\ms{O}(s^2)\,.
\]
Of course the presented conservation laws are just  special cases of Noether's Theorem.
\item For the construction of gauge theories one needs, in physicist's terminology,  a 'global symmetry of the Lagrangian'.
To achieve this,  {\em an Internal Symmetry} of the proto-Lagrangian $\vwL$ is required here: For some fixed Lie-group $\gtG\subset\dC^{c\ti c}$, the proto-Lagrangian satisfies
\[
\vwL(\mP\gU\,;\os{\mQ}\gU\,;\,\ox)=\vwL(\mP\,;\os{\mQ}\,;\,\ox)~,~~\text{for all}~\mP\in\dC^{r\ti c}\,,\os{\mQ}\in\dC^{Nr\ti c}~,~\gU\in\gtG~,~\ox\in\dR^N\,.
\]
Roughly speaking, a gauge theory for  a Lagrangian based system of PDE's is some kind of  symmetry preserving extension of the original Lagrangian {\em density} with new (dependent) 'field'-variables $\ox\mapsto\oA(\ox)=[\aA_1(\ox),\ldots,\aA_N(\ox)]$ on $\dR^N$ added,
such that the original 'quantities' $\opsi$
become subjected to the 'gauge fields' $\oA$ and viceversa. Since about a century, Weyl 1918, it is well known that, given the existence of some 'global symmetry group' $\gtG$ of $\vwL$, an extension of type \[\vwL_{\psi,A}(\ox)~=~\vwL(\opsi\,;\nabla\opsi+\opsi\!\cdot\!\oA~;\ox)+\vwG(\oA\,;\nabla\oA\,;\ox)\,,
\]
 is {\em often} possible. This extension has to exhibit
 what physicists call, a {\em 'Local Symmetry'} : The Lagrangian density remains unaltered if in $\vwL_{\psi,A}$ the quantities $\opsi$ and $\oA$ are, each in their own way,
subjected to group actions taken from
$\gtG_{\text{loc}}=\sC^\infty(\dR^N;\gtG)$, which is the {\em group} of smooth  maps $\dR^N\to\gtG$.
The added 'gauge fields' $\oA$ have to take their values in the Lie Algebra $\gtg$ of the symmetry group $\gtG$.

Summarizing, 'locally symmetric' means, symbolically,
\[
\vwL\big(\opsi U~;\nabla(\opsi U)+(\opsi U)\cdot(\oA\dhr U)\,;\ox)+\vwG(\oA\dhr U~;\nabla(\oA\dhr U)\,;\ox\,\big)~= \hspace{4cm}
\]
\[
\hspace{4cm} ~=~  \vwL\big(\opsi~;\nabla\opsi+\opsi\!\cdot\!\oA~;\ox\,\big)+\vwG\big(\oA\,;\nabla\oA\,;\ox\,\big)\,,~~\text{for all}~~U\in\gtG_{\text{loc}}\,.
\]
\item   The considerations in the underlying notes not only include  the standard hyperbolic evolution equations of pre-quantized fields. 
Wide classes of parabolic/elliptic systems turn out to have
gauge extensions as well. Note the subtle extra condition (5.14) in Thm 5.5 which is, besides internal symmetry of the proto-Lagrangian, necessary for gauge extensions.
Its necessity lies in the fact that one has to reconcile the {\em complex} vector space, in which the $\opsi$ take their values, with the {\em real} vector space $\gtg$, the Lie-Algebra. In the standard preludes to quantum field the requirement (5.14) 
is never discussed, but manifestly met with.
\item These notes do not contain functional analysis or differential geometry. The reader will find only bare elementary considerations on matrix-valued functions: The columns of the
$\ox\mapsto\opsi(\ox)\in\dC^{r\ti c}$ might describe the 'pre-quantized wave functions' of individual elementary particles, whereas the 'components' of $\ox\mapsto\oA(\ox)\in\gtg^N$, with $\gtg\subset\dC^{c\ti c}$,
might represent the pre-quantized gauge fields. For an elementary and very readable account on the differential geometrical aspects, see the contributions 3-4 in [JP].
\end{itemize}

\  \\[0cm]

\begin{center}
 {\bf CONTENTS }
\end{center}
\begin{itemize}
\item[{\bf 1.}] Foretaste: Some gauge-type calculations     \hfill{p.3}
\item[{\bf 2.}] Stationary points of complex-valued functionals  \hfill{p.6}
\item[{\bf 3.}] Free Gauge Fields  \hfill{p.13}
\item[{\bf 4.}] Noether Fluxes  \hfill{p.19}
\item[{\bf 5.}] Static/Dynamic Gauge Extensions of Lagrangians \hfill{p.26}
\item[{\bf A.}] Addendum on Free Gauge Fields \hfill{p.34}
\item[{\bf B.}] Electromagnetism \hfill{p.35}
\item[] References \hfill{p.36}
\end{itemize}
 \ \\[1mm]

\section{Foretaste: Some gauge-type calculations}
For functions $\dpsi:~\dR^N\to\dC^{r\times c}$
we consider, by way of example, the PDE 
\be
\Gamma^\mu\big(\pa_\mu\dpsi+\dpsi \cA_\mu\big)+M\dpsi=f,
\ee
with prescribed matrix valued coefficients
\[\Gamma^\mu:~\dR^N\to\dC^{r\times r},~~~\cA_\mu:~\dR^N\to\dC^{c\times c},~1\leq \mu\leq N,~~~ M:~\dR^N\to\dC^{r\times r},~
\]
and prescribed right hand side $f:~\dR^N\to\dC^{r\times c}$. All considered functions are supposed to be sufficiently smooth.
The summation convention for upper and lower indices applies.

In physics each column of $\dpsi$ may represent a 'classical-particle wave'. The $\cA_\mu$ may then represent 'gauge fields'.
\bt 
Let $\cU\,,\,\cV:~\dR^N\to\dC^{c\times c}$ and suppose them invertible with  $\cU^{-1}\,,\,\cV^{-1}:~\dR^N\to\dC^{c\times c}$.\\
The function $\hpsi=\dpsi \cU:~\dR^N\to\dC^{r\times k}$, with $\dpsi$ any solution of (1.1) is a solution of
\be
\Gamma^\mu\big(\pa_\mu\hpsi+\hpsi \hA_\mu\big)+M\hpsi=\hf,
\ee
if and only if we take  the new coefficients $\hA_\mu=\cU^{-1}\cA_\mu \cU - \cU^{-1}(\pa_\mu \cU)$ and $\hf=f\cU$.

In addition we have $\hat{\hA}_\mu=(\cU \cV)^{-1}\cA_\mu (\cU\cV) - (\cU\cV)^{-1}(\pa_\mu (\cU\cV))=\cV^{-1}\hA_\mu \cV - \cV^{-1}(\pa_\mu \cV)$.
\et
{\bf Proof:} Multiply (1.1) from the right by $\cU$ and rearrange.
\hfill{$\blacksquare$}
\\[4mm]
In the next Theorem a 'transformation property' for matrix valued functions is derived.
\bt 
Let $\cA_\mu:~\dR^N\to\dC^{c\times c}$ and $\hA_\mu=\cU^{-1}\cA_\mu \cU - \cU^{-1}(\pa_\mu \cU)$.
Define
\be
\cF_{\mu\nu}=\pa_\mu\cA_\nu-\pa_\nu\cA_\mu-\big(\cA_\mu\cA_\nu-\cA_\nu\cA_\mu\big).
\ee
Then
\be
\hF_{\mu\nu}=\pa_\mu\hA_\nu-\pa_\nu\hA_\mu-\big(\hA_\mu\hA_\nu-\hA_\nu\hA_\mu\big)=\cU^{-1}\cF_{\mu\nu}\cU.
\ee
\et
{\bf Proof:} First note that from $\pa_\mu(\cU^{-1}\cU)=\pa_\mu I=0$ it follows that $\pa_\mu(\cU^{-1})=-\cU^{-1}(\pa_\mu \cU) \cU^{-1}$.\\[1mm]
Calculate
\[ \pa_\mu\hA_\nu=\pa_\mu\big(\cU^{-1}\cA_\nu \cU - \cU^{-1}(\pa_\nu \cU)\big)=~~~~~~~~~~~~~~~~~~~~~~~~~~~~~~~~~~~~~~~~~~~~~~~~~~~~~~~~
\]
\[ = \cU^{-1}(\pa_\mu\cA_\nu) \cU-\cU^{-1}(\pa_\mu \cU) \cU^{-1}\cA_\nu \cU
+U^{-1}\cA_\nu(\pa_\mu \cU)+
\cU^{-1}(\pa_\mu \cU) \cU^{-1}(\pa_\nu \cU)-\cU^{-1}(\pa_\mu\pa_\nu \cU).
\]
and
\[ \hA_\mu\hA_\nu=\big\{\cU^{-1}\cA_\mu \cU - \cU^{-1}(\pa_\mu \cU)\big\}\big\{\cU^{-1}\cA_\nu \cU - \cU^{-1}(\pa_\nu \cU)\big\}=~~~~~~~~~~~~~~~~~~~~~~~~~~
\]
\[= \cU^{-1}\big(\cA_\mu\cA_\nu\big)\cU-\big(\cU^{-1}\cA_\mu \cU\big)\big(\cU^{-1}(\pa_\nu \cU)\big)-\big(\cU^{-1}(\pa_\mu \cU)\big)\big(\cU^{-1}\cA_\nu \cU\big)
+\big(\cU^{-1}(\pa_\mu \cU)\big)\big(\cU^{-1}(\pa_\nu \cU)\big).
\]
Interchange the indices for two more terms and add according to (1.4). All rubbish terms cancel out.\hfill{$\blacksquare$}
\\[7mm]
We now look for sesqui-linear conservation laws which hold for suitable classes of $\cA_\mu$
\bc
 $K:~\dR^N\to\dC^{r\times r},~$ is such that \\{\bf i}: $K\Gamma^\mu=(K\Gamma^\mu)^\dagger,~~~$ {\bf ii}: $\pa_\mu(K\Gamma^\mu)=0,~~~$
{\bf iii}: $KM+M^\dagger K^\dagger=0$.
\ec
Here, the dagger $\dagger$ denotes 'Hermitean transposition'.

Note that in the important special case that  $\Gamma^\mu=(\Gamma^\mu)^\dagger$, $\Gamma^\mu$ is constant and $M=-M^\dagger$, the condition is satisfied by $K=I$, the identity matrix. In the case of the Dirac equation one could take $K=\Gamma^0$. Cf. [M], Messiah II pp. 890-899.
\footnote{In the non-covariant form, i.e. the original form, of Dirac's equation one has $\Ga^0=I, \Ga^\ka=\ga^0\ga^\ka\,,\,1\leq\ka\leq3\,$,
where the $\ga^\mu\,,\,0\leq\mu\leq3$ are Dirac-Clifford matrices, which make the Dirac equation covariant proof.}

\bt
Let $K:~\dR^N\to\dC^{r\times r}~$ satisfy Condition 1.3.

Fix some $J\in \dC^{c\times c}$.\\
Let $\cA_\mu:~\dR^N\to\dC^{c\times c}$ satisfy $\cA_\mu^\da J+ J\cA_\mu=0,~1\leq \mu\leq N$. \\
Let  $U:~\dR^N\to\dC^{c\times c}~$ satisfy $~U^\da(\osx) J U(\osx)=J\,,~~\osx\in\dR^N$.\\
{\bf a.} For any solution $\dpsi$ of (1.1) with $f=0$, there is the conservation law
\be 
\sum_{\mu=1}^N\,\pa_\mu\Tr\big( J^{-1}[ \dpsi^\dagger K\Gamma^\mu\dpsi]\big)=0\,.
\ee
{\bf b.} This conservation law is  a {\em gauge invariant local conservation law}. \\
That means $\Tr\big( J^{-1}[ \hpsi^\dagger K\Gamma^\mu\hpsi]\big)=\Tr\big( J^{-1}[ \dpsi^\dagger K\Gamma^\mu\dpsi]\big)~,\,1\leq\mu\leq N$.
\et
{\bf Proof}\\
{\bf a.} Take $f=0$ in (1.1)and multiply  from the left with $\dpsi^\dagger K$:
\be
\dpsi^\dagger K\Gamma^\mu\big(\pa_\mu\dpsi\big)+\dpsi^\dagger K\Gamma^\mu\dpsi \cA_\mu +\dpsi^\dagger K M\dpsi~ =~ 0.
\ee
The Hermitean transpose reads
\be
\big(\pa_\mu\dpsi\big)^\dagger(K\Gamma^\mu)^\dagger\dpsi+\cA_\mu^\dagger\dpsi^\dagger(K\Gamma^\mu)^\dagger\dpsi +\dpsi^\dagger M^\dagger K^\dagger\dpsi =0.
\ee
Multiply (1.6) from the right with $J^{-1}$ and (1.7) from the left with $J^{-1}$.
Add those two identities and take the trace. Use Condition 1.3 and the properties $\Tr(AB)=\Tr(BA)$,  $\Tr(A+B)=\Tr(A)+\Tr(B)$ and $\pa_\mu\Tr(A)=\Tr(\pa_\mu A)$.
The sum of the 1st terms of (1.6), (1.7) result in
\[
\Tr \big\{J^{-1}\big[ \dpsi^\dagger (K\Gamma^\mu)\pa_\mu\dpsi+(\pa_\mu\dpsi)^\dagger (K\Gamma^\mu)^\dagger\dpsi\big]\big\}=~~~~~~~~~~~~~~~~~~~~~~~~~~~~~~~~~~~~~~~~~~~~~~~~~~~
\]
\[=\pa_\mu\Tr\big\{J^{-1}\dpsi^\dagger (K\Gamma^\mu)\dpsi\big\}-\Tr\big\{J^{-1}\dpsi^\dagger\pa_\mu(K\Gamma^\mu)\dpsi\big\}=\pa_\mu\Tr\big\{J^{-1}\dpsi^\dagger (K\Gamma^\mu)\dpsi\big\}.
\]
The sum of the 2nd terms of (1.6), (1.7) is
\[
\Tr\big\{\dpsi^\dagger K\Gamma^\mu\dpsi\big( \cA_\mu J^{-1}+J^{-1}\cA_\mu^\dagger\big)\big\}=0.
\]
The sum of the 3rd terms of (1.6), (1.7)
\[\Tr\big\{ J^{-1}\dpsi^\dagger( KM+M^\da K^\da)\dpsi\big\}=0.
\]
Thus, we find (1.5)
\\ {\bf b.} By putting hats on $\dpsi$ and $\cA_\mu$ our considerations can be rephrased for PDE (1.2). Remind that from $U^\da JU=J$ it follows that
$J^{-1}U^\da=U^{-1}J^{-1}$. Finally
\[ 
\Tr\big( J^{-1}U^\da[ \dpsi^\da K\Gamma^\mu\dpsi]U\big)=\Tr\big( U^{-1}J^{-1}[ \dpsi^\dagger K\Gamma^\mu\dpsi]U\big)=
\Tr\big(J^{-1}[ \dpsi^\dagger K\Gamma^\mu\dpsi]\big).
\]
\hfill{$\blacksquare$}

\section{Stationary points of complex-valued functionals}

In this section we pay some attention to the Euler Lagrange field equations in the {\em complex field} case.  Most  physics textbooks start,
in a rather verbose way,
with 18th century variational  rituals.  However most of them become suddenly  very vague, or fall completely silent, when  state functions involving {\em complex variables} come into play! In order to get some feeling for  such Lagrangians, we first mention  a finite dimensional toy result.
\bt
Let
\[ f:~ \dC^n\times\dC^n~\ni~(\osz;\osw)\mapsto f(\osz,\osw)~\in~\dC~~~~
\]
be an {\em analytic} function of $2n$ complex variables with the special property $f(\osz,\osz^\star)\in \dR$, for all $\osz\in\dC^n$.
Here $\osz=\osx + \ri\osy$, $~\osz^\star=\osx - \ri\osy$.

{\bf a.} Consider the function
\[
\dR^n\times\dR^n~\ni~(\osx;\osy)\mapsto g(\osx,\,\osy)=f(\osz,\,\osz^\star)=f(\osx+\ri\osy,\,\osx-\ri\osy)~\in~\dR.
\]
The relations between the (real) partial derivatives of $g$ at $(\osx,\osy)$ and the (complex) partial derivatives of $f$ at $(\osz,\osz^\star)$ are
\be\displaystyle
\begin{array}{lcl}
\displaystyle\frac{\pa g}{\pa\osx}(\osx,\,\osy)=\frac{\pa f}{\pa\osz}(z,\,z^\star)+\frac{\pa f}{\pa\osw}(z,\,z^\star) &~~~~~
& \displaystyle\frac{\pa f}{\pa\osz}(\osz,\,\osz^\star)=\frac12\big(\frac{\pa g}{\pa\osx}(\osx,\,\osy)-\ri\frac{\pa g}{\pa\osy}(\osx,\,\osy)\big)
\\[5mm]
\displaystyle\frac{\pa g}{\pa\osy}(\osx,\,\osy)=\ri\frac{\pa f}{\pa\osz}(z,\,z^\star)-\ri\frac{\pa f}{\pa\osw}(z,\,z^\star) &~~~~
& \displaystyle\frac{\pa f}{\pa\osw}(z,\,z^\star)=\frac12\big(\frac{\pa g}{\pa\osx}(\osx,\,\osy)+\ri\frac{\pa g}{\pa\osy}(\osx,\,\osy)\big)\,          \\
\end{array}
\ee
\[
\frac{\pa f}{\pa \osw}(\osz,\osz^\star)=\bs{\frac{\pa f}{\pa \osz}(\osz,\osz^\star)}
\]
{\bf b.} For $g$ to have a stationary point at $(\osa\,;\,\osb)\in\dR^n\times\dR^n$ {\bf each one} of the following three  conditions is necessary and sufficient
\be
\begin{array}{l}
\bullet~~~\displaystyle\frac{\pa g}{\pa\osx}(\osa,\,\osb)=\frac{\pa g}{\pa\osy}(\osa,\,\osb)=\onul\,,\\[4mm]
\bullet~~~\displaystyle\frac{\pa f}{\pa \osz}(\osa+\ri\osb,\,\osa-\ri \osb)=\onul\,, \\[4mm]
\bullet~~~\displaystyle\frac{\pa f}{\pa \osw}(\osa+\ri \osb,\,\osa-\ri \osb)=\,"\,\frac{\pa f}{\pa \osz^\star}(\osa+\ri \osb,\,\osa-\ri \osb)\,"\,=\onul\,.\\
\end{array}
\ee
{\bf c.} If the special property $f(\osx+\ri\osy,\,\osx-\ri\osy)\in\dR$ is relaxed to $\phi(f(\osx+\ri\osy,\,\osx-\ri\osy))\in\dR$ for some non-constant analytic $\phi:\dC\to\dC$,
then the 'stationary point result' {\bf b.} still holds.
\et
{\bf Proof:} Straightforward calculation \hfill{$\blacksquare$}\\[2mm]
In Theorem 2.4 an $\infty$-dimensional generalisation of this result is presented.\\[3mm]

{\bf A special bookkeeping} 

In the sequel, for the above variable $\osz$,  usually  a matrix $\mZ\in\dC^{r\ti c}$ will be taken. In order to explain our bookkeeping and also for some special properties, we now consider an analytic function of 2 matrix variables
\be
\sF:~~\dC^{r\times c}\times\dC^{c\times r}~\to~\dC~~:~~(\mZ\,;\mW)~\mapsto~\sF(\mZ\,,\mW)\,.
\ee
Because of Hartog's Theorem, see [H] Thm 2.2.8, it is enough to assume
  analyticity with respect to each entry of each matrix separately.

 The (complex!) partial derivatives of $\sF$ are gathered in matrices,
 \[
 (\mZ;\mW)~\mapsto~\sF^{(\dg{1})}(\mZ,\mW)~\in~\dC^{c\times r}\,, ~~~~~~(\mZ;\mW)~\mapsto~\sF^{(\dg{2})}(\mZ,\mW)~\in~\dC^{r\times c}\,,
 \]

 with
 \be
 \big[\sF^{(\dg{1})}\big]_{ij}=\big[\dfrac{\pa\sF}{\pa\mZ}\big]_{\mbox{}_{ij}}=\dfrac{\pa\sF}{\pa \mZ_{ji}}~,~~~~~~~~~~
 \big[\sF^{(\dg{2})}\big]_{k\ell}=\big[\dfrac{\pa\sF}{\pa\mW}\big]_{\mbox{}_{k\ell}}=\dfrac{\pa\sF}{\pa \mW_{\ell k}}.
 \ee

 In our notation the $\dC$-linearization of $\sF$ at $(\mZ,\mW)$, for $\vep\in\dC\,,|\vep|$ small, reads

 \be
 \sF(\mZ+\vep\mH,\mW+\vep\mK)=\sF(\mZ,\mW)+\vep\Tr\big\{[\sF^{(\dg{1})}]\mH\}+\vep\Tr\big\{[\sF^{(\dg{2})}]\mK\big\}+\mathscr{O}(|\vep|^2).
 \ee
 {\bf Notation:} Sometimes, in order to avoid excessive use of brackets, it is convenient to write \\ $\Tr\big\{\sF^{(\dg{1})}:\mH\}$ instead of $\Tr\big\{[\sF^{(\dg{1})}]\mH\}$.\\
 Also, without warning, in proofs sometimes Einstein's summation convention for repeated upper and lower indices will be used.\\[0.5mm]

 Next split $\mZ$ in real and imaginary parts $\mZ=\mX+\ri\mY$ and introduce the function
 \be
 \ssF:~~\dR^{r\times c}\times\dR^{r\times c}~\to~\dC~~:~~(\mX;\mY)~\mapsto~\ssF(\mX\,,\mY)=\sF(\mZ,\mZ^\dagger)=\sF(\mX+\ri\mY\,,\mX^\top-\ri\mY^\top).
 \ee
 The $\dR$-linearization of $\ssF$ at $(\mX\,,\mY)$ for $\vep\in\dR\,,|\vep|$ small, can now be written
 \be
 \ssF(\mX+\vep\mA\,,\mY+\vep\mB)=\ssF(\mX\,,\mY)+\vep\Trc{\dfrac{\pa\ssF}{\pa\mX}\mA+\vep\Trc{\dfrac{\pa\ssF}{\pa\mY}\mB}} +\mathscr{O}(\vep^2),
 \ee
with
\be
\begin{array}{l}
\Trc{\dfrac{\pa\ssF}{\pa\mX}\mA}=\Trc{[\sF^{(\dg{1})}]\mA}+\Trc{[\sF^{(\dg{2})}]\mA^\top}=\Trc{\big(\,[\sF^{(\dg{1})}]+[\sF^{(\dg{2})}]^\top\big)\mA},\\[3mm]
\Trc{\dfrac{\pa\ssF}{\pa\mY}\mB}=\Trc{\ri[\sF^{(\dg{1})}]\mB}+\Trc{-\ri[\sF^{(\dg{2})}]\mB^\top}=\Trc{\ri\big(\,[\sF^{(\dg{1})}]-[\sF^{(\dg{2})}]^\top\big)\mB},
\end{array}
\ee
where the matrices $\mX,\mY,\mA,\mB$ are all real. The (complex) derivatives  $\sF^{(\dg{1})}\,,\sF^{(\dg{2})}$ are taken at $(\mZ,\mZ^\dagger)$. In the usual (somewhat confusing) notation, this corresponds to
\be
\dfrac{\pa\ssF}{\pa\mX}=\dfrac{\pa\sF}{\pa\mX}=\dfrac{\pa\sF}{\pa\mZ}+\big[\dfrac{\pa\sF}{\pa\mZ^\dagger}\big]^\top~,~~~~
\dfrac{\pa\ssF}{\pa\mY}=\dfrac{\pa\sF}{\pa\mY}=\ri\dfrac{\pa\sF}{\pa\mZ}-\ri\big[\dfrac{\pa\sF}{\pa\mZ^\dagger}\big]^\top\,,
\ee
and, similarly sloppy,
\be
\dfrac{\pa\sF}{\pa\mZ}=\frac12\big(\dfrac{\pa\sF}{\pa\mX}-\ri\dfrac{\pa\sF}{\pa\mY}\big)~,~~~~
\big[\dfrac{\pa\sF}{\pa\mZ^\dagger}\big]^\top=\frac12\big(\dfrac{\pa\sF}{\pa\mX}+\ri\dfrac{\pa\sF}{\pa\mY}\big).
\ee
\ \\[4mm]
If it happens that $\mZ\mapsto\sF(\mZ,\mZ^\dagger)$ is $\dR$-valued, the results of Theorem (2.1) can be rephrased.
\bt
Let, as in (2.3),
\[ \sF:~~\dC^{r\times c}\times\dC^{c\times r}~~\ni~~(\mZ;\mW)~\mapsto~\sF(\mZ,\mW)~\in~\dC~.
\]
be  {\em analytic}. Suppose $\sF(\mZ,\mZ^\dagger)\in\dR$, for all $\mZ\in\dC^{r\times c}$. Write $\mZ=\mX+\ri\mY$. Denote
\[
 \ssF:~~\dR^{r\times c}\times\dR^{r\times c}~\to~\dR~~:~~(\mX;\mY)~\mapsto~\ssF(\mX\,,\mY)=\sF(\mZ,\mZ^\dagger)=\sF(\mX+\ri\mY\,,\mX^\top-\ri\mY^\top)\,,
\]
$\bullet$ We have
\be \sF^{(\dg{1})}(\mZ,\mZ^\dagger)\,=\,[\sF^{(\dg{2})}(\mZ,\mZ^\dagger)]^\dagger.
\ee

Further, for the function $\ssF$ to have a stationary point at $(\mA\,;\,\mB)\in\dR^{r\times c}\times\dR^{r\times c}$ {\bf each one} of the following three  conditions is necessary and sufficient
\be
\begin{array}{l}
\bullet~~~\dfrac{\pa\ssF}{\pa\mX}(\mA,\mB)=\dfrac{\pa\ssF}{\pa\mY}(\mA,\mB)=0\\[4mm]
\bullet~~~\displaystyle\sF^{(\dg{1})}(\mA+\ri\mB,\,\mA^\top-\ri \mB^\top)=\,"\,\frac{\pa \sF}{\pa \mZ}(\mA+\ri\mB,\,\mA^\top-\ri \mB^\top)\,"=0\\[4mm]
\bullet~~~\displaystyle\sF^{(\dg{2})}(\mA+\ri\mB,\,\mA^\top-\ri \mB^\top)=\,"\,\frac{\pa \sF}{\pa \mZ^\dagger}(\mA+\ri\mB,\,\mA^\top-\ri \mB^\top)\,"\,=0.
\end{array}
\ee
\et

{\bf Proof:}
Is mostly a reformulation of the preceding theorem.  It follows directly from (2.9)-(2.10). \eindebewijs
\ \\[5mm]

In order to build the concept of { \bf Lagrangian density} we need an analytic function, named { \bf proto-Lagrangian},
\be\begin{array}{c}
\sL:~~\dC^{r\times c}\times\dC^{c\times r}\times\dC^{Nr\times c}\times\dC^{c\times Nr}\times\dR^N~~~\to~~~\dC,\\[1mm]

(\mP;\mQ^\top;\os{\mR}\,;\os{\mS}^\top;\,\osx)~~\mapsto~~\sL(\mP;\mQ^\top;\os{\mR}\,;\os{\mS}^\top;\,\osx)\,,\\
\end{array}
\ee
where
\[
\begin{array}{llc}
\mP\in\dC^{r\times c}~,~&   \os{\mR}=\text{col}\big[\,\mR_1\,,\ldots,\mR_N\,\big]~,~&\mR_\mu\in\dC^{r\times c}~,~1\leq \mu\leq N\,,\\[1mm]

\mQ^\top\in\dC^{c\times r}~,~& \os{\mS}^\top=\text{row}\big[\,\mS^\top_1\,,\ldots,\mS_N^\top\,\big]~,~&\mS_\mu^\top\in\dC^{c\times r}~,~1\leq \mu\leq N\,.\\
\end{array}
\]
Instead of (2.13) it will be convenient sometimes to denote the proto Lagrangian by
\[
\sL(\mP;\mQ^\top;\ldots,\os{\mR}_\mu,\ldots\,;\ldots,\os{\mS}^\top_\mu,\ldots;\,\osx).
\]
It will be required that $\sL(\mO;\mO^\top;\os{\mO}\,;\os{\mO}^\top;\,\osx)=0$.
\\[1mm]
The (complex) partial derivatives of $\sL$, cf. (2.4)-(2.5), with respect to its $2N+2$ matrix arguments are denoted, respectively,
\[
\sL^{\dg{(o)}}\,,\,\sL^{\dg{(o\star)}}\,,\,\sL^{\dg{(1)}}\,,\ldots\,,\sL^{\dg{(N)}}\,,\,\sL^{\dg{(1\star)}}\,,\ldots\,,\sL^{\dg{(N\star)}}\,.
\]
The (real) partial derivatives of $\sL$, with respect to the vector variable $\osx$ is denoted $~\sL^{\dg{(\nabla)}}~$.
For any given
 matrix-valued function $\dpsi:~\dR^N\to\dC^{r\times c}$,  we define a {\em Lagrangian density} $\sL_\psi:~\dR^N\to\dC$, by substitution of  $\dpsi$, its 1st derivatives $\pa_\mu\dpsi=\dpsi_{,\,\mu}\,,~1\leq \mu\leq N$, and the hermitean transposed of all those, in $\sL$:

\be
\osx~~\mapsto~~\sL_\psi(\osx)=\sL(\dpsi(\osx);\dpsi^\dagger(\osx);\nabla\dpsi(\osx)\,;\nabla\dpsi^\dagger(\osx)\,;\osx\,),
\ee
where
\[
\nabla\dpsi(\osx)=\text{col}\big[\pa_1\dpsi(\osx)\,,\dots,\pa_N\dpsi(\osx)\big]\,\in\,\dC^{Nr\times c}\,,
\]
\[
\nabla\dpsi^\dagger(\osx)=\text{row}\big[\pa_1\dpsi^\dagger(\osx)\,,\dots,\pa_N\dpsi^\dagger(\osx)\big]\,\in\,\dC^{c\times Nr}\,.
\]
Also the matrix-valued functions
\[
\osx~~\mapsto~~[\sL^{\dg{(\mu)}}_\psi](\osx)=[\sL^{\dg{(\mu)}}](\dpsi(\osx);\dpsi^\dagger(\osx);\nabla\dpsi(\osx)\,;\nabla\dpsi^\dagger(\osx)\,;\osx\,)\,\in\,\dC^{c\ti r},
\]
similarly $\osx\mapsto[\sL^{\dg{(\mu\star)}}_\psi]\,\in\dC^{r\ti c}$, and $\osx\mapsto\sL^{\dg{(\nabla)}}_\psi\,\in\dR^N\,$, will be used.
\  \\[2mm]
On a suitable  space of functions $\dpsi:\dR^N\to\dC^{r\times c}$, it often makes sense to define the {\em \bf Lagrangian functional}
 \be\dpsi~~\mapsto~~\cL(\dpsi\,,\dpsi^\dagger)=\int_{\dR^N}\,\sL(\dpsi(\osx);\dpsi^\dagger(\osx);\nabla\dpsi(\osx)\,;
 \nabla\dpsi^\dagger(\osx)\,;\osx\,)\rd\osx\,\in\,\dC.
\ee
\begin{remarkn}
The Lagrangian functional $\cL$ remains the same if we replace $\sL$ by
\[ \sL(\dpsi;\dpsi^\dagger;\nabla\dpsi\,;\nabla\dpsi^\dagger;\,\osx)+ \pa_\mu w^\mu(\dpsi,\dpsi^\dagger,\,\osx),
\]
with $w^\mu$ a vectorfield which vanishes sufficiently rapidly at infinity.\\[1mm]
Therefore the functional $\dpsi~~\mapsto~~\cL(\dpsi\,,\dpsi^\dagger)$ is $\dR$-valued if
\[
\overline{\sL(\dpsi;\dpsi^\dagger;\nabla\dpsi\,;\nabla\dpsi^\dagger;\,\osx)}\,-\,\sL(\dpsi;\dpsi^\dagger;\nabla\dpsi\,;\nabla\dpsi^\dagger;\,\osx)\,=\,\pa_\mu W^\mu(\dpsi,\dpsi^\dagger,\,\osx)\,,
\]
i.e. the divergence of a vector field.

{\bf Note} that $\cL$ may be $\dR$-valued while $\sL_\psi$ is not\,!!
\end{remarkn}
If we split    $\dpsi$  into real and imaginary parts: $\dpsi=\dpsi_\Re + \ri\dpsi_\Im$ and $\dpsi_{,\mu}=\dpsi_{\Re,\mu} + \ri\dpsi_{\Im,\mu}~$,

the $\dR$-directional derivatives with respect to $\dpsi_\Re$ and $\dpsi_\Im$ of the Lagrangian functional $\cL$ are explained by
\[
\begin{array}{l}
\big\langle\,\Der_{\dpsi_\Re}\cL\,,\,\dk{A}\,\big\rangle  = \dfrac{\rd}{\rd \vep}\cL(\dpsi+\vep \dk{A}\,,\,\dpsi^\dagger+\vep \dk{A}^\top)\Big|_{\vep=0}=\\[3mm]
\  \hspace{-9mm}=\dfrac{\rd}{\rd \vep}\dint_{\!\!\dR^N}\,\sL(\dpsi(\osx)+\vep\dk{A}(\osx);\dpsi^\dagger(\osx)+\vep\dk{A}^\top(\osx);\nabla\big(\dpsi(\osx)+\vep\dk{A}(\osx)\big)
\,;\nabla\big(\dpsi^\dagger(\osx)+\vep\dk{A}^\top(\osx)\big)\,;\osx)\rd\osx\Big|_{\vep=0},\\[3mm]
\hspace{40mm}\mbox{with}~~\dk{A}~:~\dR^N\to\dR^{r\times c}\,,\,\text{and}~\vep\in\dR\,,\,|\vep|~~\text{small}.\\[3mm]

\big\langle\,\Der_{\dpsi_\Im}\cL\,,\,\dk{B}\,\big\rangle  = \dfrac{\rd}{\rd \vep}\cL(\dpsi+\vep \ri\dk{B}\,,\,\dpsi^\dagger-\vep \ri\dk{B}^\top)\Big|_{\vep=0}=\\[3mm]
\ \hspace{-9mm}=\dfrac{\rd}{\rd \vep}\dint_{\!\!\dR^N}\,\sL(\dpsi(\osx)+\vep \ri\dk{B}(\osx);\dpsi^\dagger(\osx)-\vep \ri\dk{B}^\top(\osx);\nabla\big(\dpsi(\osx)+\vep \ri\dk{B}(\osx)\big)
\,;\nabla\big(\dpsi^\dagger(\osx)-\vep \ri\dk{B}^\top(\osx)\big)\,;\osx)\rd\osx\Big|_{\vep=0},\\[3mm]
\hspace{40mm}\mbox{with}~~\dk{B}~:~\dR^N\to\dR^{r\times c}\,,\,\text{and}~\vep\in\dR\,,\,|\vep|~~\text{small}.\\[1mm]
\end{array}
\]

When calculating the $\dC$-directional derivatives  $\Der_{\dpsi}\cL\,,\Der_{\dpsi^\dagger}\cL\,$, the variables $\dpsi\,,\dpsi^\dagger$ are considered to be independent. These derivatives are supposed to be elements in the {\em (complex) linear} dual of $\dL_2(\dR^N;\dC^{r\times c})$. They are explained by

\[
\begin{array}{l}
\big\langle\,\Der_{\dpsi}\cL\,,\,\dk{H}\,\big\rangle  = \dfrac{\rd}{\rd \vep}\cL(\dpsi+\vep \dk{H}\,,\,\dpsi^\dagger)\Big|_{\vep=0}=\\[3mm]
~~~~~~~~=\dfrac{\rd}{\rd \vep}\dint_{\!\!\dR^N}\,\sL(\dpsi(\osx)+\vep\dk{H}(\osx);\dpsi^\dagger(\osx);\nabla\big(\dpsi(\osx)+\vep\dk{H}(\osx)\big)
\,;\nabla\dpsi^\dagger\,;\,\osx\,)\rd\osx\Big|_{\vep=0},\\[3mm]
\hspace{40mm}\mbox{with}~~\dk{H}~:~\dR^N\to\dC^{r\times c}\,,\,\text{and}~\vep\in\dC\,,\,|\vep|~~\text{small}.\\[3mm]

\big\langle\,\Der_{\dpsi^\dagger}\cL\,,\,\dk{K}\,\big\rangle  = \dfrac{\rd}{\rd \vep}\cL(\dpsi\,,\,\dpsi^\dagger+\vep \dk{K})\Big|_{\vep=0}=\\[3mm]
~~~~~~~~=\dfrac{\rd}{\rd \vep}\dint_{\!\!\dR^N}\,\sL(\dpsi(\osx);\dpsi^\dagger(\osx)+\vep\dk{K}(\osx);\nabla\dpsi(\osx)
\,;\nabla(\dpsi^\dagger(\osx)+\vep\dk{K}(\osx))\,;\osx\,)\rd\osx\Big|_{\vep=0},\\[3mm]
\hspace{40mm}\mbox{with}~~\dk{K}~:~\dR^N\to\dC^{c\times r}\,,\,\text{and}~\vep\in\dC\,,\,|\vep|~~\text{small}.\\[3mm]
\end{array}
\]
For $\dk{H}\,,\dk{K}\,,\dk{A}\,,\dk{B}$ vanishing sufficiently rapidly at $\infty$ a partial integration leads to the standard Euler-Lagrange expressions for the functional derivatives of $\cL$.

\bt
Assume that  $\cL$ is $\dR$-valued. (Cf. Remark 2.3). If $\dpsi$ satisfies any one of the following three Lagrangian systems
\be
\begin{array}{cc}

\begin{array}{l}
\displaystyle{
\Der_{\dpsi}\cL=[\sL^{\dg{(o)}}_\psi]-\sum_{\mu=1}^N\frac{\pa}{\pa x^\mu}[\sL^{\dg{(\mu)}}_\psi]=0~,
}
\\[5mm]
\displaystyle{
\Der_{\dpsi^\dagger}\cL=[\sL^{\dg{(o\star)}}_\psi]-\sum_{\mu=1}^N \frac{\pa}{\pa x^\mu}[\sL^{\dg{(\mu\star)}}_\psi]=0~,}\\
  \end{array}
&
~~~\left\{
\begin{array}{l}
\displaystyle{
\Der_{\dpsi_{\Re}}\cL=\frac{\pa\sL~}{\pa\dpsi_\Re}-\sum_{\mu=1}^N\frac{\pa}{\pa x^\mu}\frac{\pa\sL~}{\pa \dpsi_{\Re,\mu}}=0~,}
\\[4mm]
\displaystyle {
\Der_{\dpsi_{\Im}}\cL=\frac{\pa\sL~}{\pa\dpsi_\Im}-\sum_{\mu=1}^N\frac{\pa}{\pa x^\mu}\frac{\pa\sL~}{\pa \dpsi_{\Im,\mu}}=0~.
}
\end{array}
\right.
\end{array}~~,
\ee
with $\sL=\sL(\dpsi(\osx);\dpsi^\dagger(\osx);\nabla\dpsi(\osx)\,;\nabla\dpsi^\dagger(\osx)\,;\osx\,)\,$, then it also satisfies the other two. 
\et
{\bf Proof:} With the notation (2.8)-(2.10) we obtain
\be
\frac{\pa\sL~}{\pa\dpsi_\Re}=\sL^{\dg{(o)}}+[\sL^{\dg{(o\star)}}]^\top~~,~~~~~
\frac{\pa\sL~}{\pa\dpsi_\Im}=\ri\sL^{\dg{(o)}}-\ri[\sL^{\dg{(o\star)}}]^\top~~,  
\ee
and, the other way round,
\be
\big[\sL^{\dg{(o\star)}}\big]^\top=
\frac12\big(\frac{\pa\sL~}{\pa\dpsi_\Re}+\ri\frac{\pa\sL~}{\pa\dpsi_\Im}\big)~~,~~~
\sL^{\dg{(o)}}=
\frac12\big(\frac{\pa\sL~}{\pa\dpsi_\Re}-\ri\frac{\pa\sL~}{\pa\dpsi_\Im}\big),
\ee
and similar expressions with  $\dg{(o)}\,,\dg{(o\star)}$ replaced by $\dg{(\mu)}\,,\dg{(\mu\star)}$ and $\dpsi\,,\dpsi_\Re\,,\dpsi_\Im$ replaced by $\dpsi_{,\mu}\,,\dpsi_{\Re,\mu}\,,\dpsi_{\Im,\mu}$.
Then
\[
\begin{array}{ccc}
\begin{array}{ccc}
\Der_{\dpsi}\cL &=&\frac12\big(\Der_{\dpsi_{\Re}}\cL-\ri\Der_{\dpsi_{\Im}}\cL\big)\\[2mm]
\big[\Der_{\dpsi^\dagger}\cL\big]^\top &=&\frac12\big(\Der_{\dpsi_{\Re}}\cL+\ri\Der_{\dpsi_{\Im}}\cL\big)\\
\end{array}
&~~~~&

\begin{array}{ccc}
\Der_{\dpsi_\Re}\cL &=&\Der_{\dpsi}\cL+\big[\Der_{\dpsi^\dagger}\cL\big]^\top\\[2mm]
\big[\Der_{\dpsi_\Im}\cL\big]^\top&=&\ri\Der_{\dpsi}\cL-\ri\big[\Der_{\dpsi^\dagger}\cL\big]^\top\\
\end{array}
\end{array}.
\]
If we take into account that the entries of the matrix valued functions $\Der_{\dpsi_\Re}\cL \,$ and $\,\Der_{\dpsi_\Im}\cL $ are $\dR$-valued, we find
\be
\big[\Der_{\dpsi^\dagger}\cL\big]^\dagger\,=\,\big[\Der_{\dpsi}\cL\big]\,,
\ee
from which the theorem easily follows. \hfill{$\blacksquare$}
\ \\

\begin{examps} { \ \bf (Matter Fields)} \\
{\bf a)} Let $\Gamma^\mu$ and $M$ be constant complex matrices with $\Gamma^{\mu\dagger}=\Gamma^\mu$ and $M=-M^\dagger$. Then the Lagrangian density
\be
\sL_\psi=\ri\,\Tr\big\{\dpsi^\dagger \Gamma^\mu\pa_\mu\dpsi+\dpsi^\dagger M\dpsi \big\},
\ee
for $\dpsi:~\dR^N\to\dC^{r\times c}$, satisfies the condition of Theorem (2.4) and leads to (1.1) with $\cA=0$.\\[2.5mm]
{\bf b)}
Let $\Gamma_\mu\,,1\leq\mu\leq N :~\dR^N\to\dC^{r\times r}$. Let $\cA_\mu\,,1\leq\mu\leq N :~\dR^N\to\dC^{c\times c}$.\\ Let $M :~\dR^N\to\dC^{r\times r}$.\\
Suppose both the existence of $K:~\dR^N\to\dC^{r\times r},~$ having inverse $K^{-1}(\osx)$,  for all $\osx\in\dR^N$, and an invertible $J\in \dC^{c\times c}$ with $J^\dagger=J$, such that:\\[2mm]
$(K\Gamma^\mu)^\dagger=K\Gamma^\mu,~,~1\leq \mu\leq N,~$ $~\cA_\mu^\dagger(\osx) J+ J\cA_\mu(\osx)=0,~1\leq \mu\leq N,~\osx\in\dR^N,~$ \\and $KM+M^\dagger K^\dagger-\pa_\mu\big(K\Gamma^\mu\big)=0$.\\[0.5mm]

Then the Lagrangian density
\be
\sL_\psi=\ri\,\Tr\big\{\dpsi^\dagger K(\Gamma^\mu\pa_\mu\dpsi) J^{-1}+ \dpsi^\dagger K(\Gamma^\mu\dpsi\cA_\mu) J^{-1}+\dpsi^\dagger KM\dpsi J^{-1} \big\},
\ee
for $\dpsi:~\dR^N\to\dC^{r\times c}$ satisfies $\sL-\bs{\sL}=\pa_\mu w$ and hence the condition of Theorem (2.4).

It leads to the 'matter-field equation'
\be
\Gamma^\mu\pa_\mu\dpsi + \Gamma^\mu\dpsi\cA_\mu +M\dpsi =0
\ee
{\bf Indeed.} Taking suitable combinations  we find respectively
\[\begin{array}{rcl}
\Tr\big\{ \dpsi^\dagger K\Gamma^\mu(\pa_\mu\dpsi) J^{-1}+J^{-1}(\pa_\mu\dpsi)^\dagger (K\Gamma^\mu)^\dagger\dpsi)\big\} &=&
  \Tr\big\{J^{-1}\pa_\mu[ \dpsi^\dagger K\Gamma^\mu\dpsi)]\big\}+\\[1mm]

  & &~~~~~~~\Tr\big\{J^{-1}[ \dpsi^\dagger \pa_\mu(K\Gamma^\mu)\dpsi)]\big\},\\[3mm]

\Tr\big\{\dpsi^\dagger K(\Gamma^\mu\dpsi\cA_\mu) J^{-1} +J^{-1}\cA_\mu^\dagger\dpsi^\dagger(K\Gamma^\mu)^\dagger\dpsi\big\}&=&
\Tr\big\{\big[\cA_\mu J^{-1}+J^{-1}\cA_\mu^\dagger\big]\dpsi^\dagger( K\Gamma^\mu)\dpsi\big\}=0,\\[3mm]

 \Tr\big\{\dpsi^\dagger KM\dpsi J^{-1} +J^{-1}\dpsi^\dagger M^\dagger K^\dagger\dpsi\big\}&=&
   \Tr\big\{J^{-1}\dpsi^\dagger KM\dpsi  +J^{-1}\dpsi^\dagger M^\dagger K^\dagger\dpsi\big\}= \\[1mm]
   &=&  \Tr\big\{J^{-1}\big[\dpsi^\dagger (KM+ M^\dagger K^\dagger)\dpsi\big]\big\}.\\
\end{array}
\]
Ultimately we find
\be
\sL_\psi - \overline{\sL_\psi}=\pa_\mu\Tr\big\{J^{-1}[ \dpsi^\dagger K\Gamma^\mu\dpsi)]\big\}=\pa_\mu\Tr\big\{[ \dpsi^\dagger K\Gamma^\mu\dpsi)]J^{-1}\big\}.
\ee
The Euler-Lagrange equations are
\be
K\big(\Gamma^\mu\pa_\mu\dpsi+  \Gamma^\mu\dpsi\cA_\mu+ M\dpsi\big) J^{-1}=0,
\ee
from which $K$ and $J^{-1}$ can be cancelled.\\[1mm]

{\bf c)}
The Lagrangian density
\be
\sL_\psi=\Tr\big\{[\pa_\mu\dpsi]^\dagger\Theta^{\mu\nu}[\pa_\nu\dpsi]+\dpsi^\dagger R\dpsi\big\}\,,
\ee
with $\Theta^{\mu\nu}, R\,:\dR^N \to\dC^{r\ti r}$ and $[\Theta^{\mu\nu}]^\dagger=\Theta^{\nu\mu}\,,R^\dagger=R$,  is $\dR$-valued. It leads to the 2nd order equation
\be
\sum_{\mu,\nu}\paf{}{x^\mu}\,\Theta^{\mu\nu}\paf{}{x^\nu}\dpsi\,-R\,\dpsi=0\,.
\ee

{\bf d.} The Lagrangian density for functions $\dpsi=\col[\begin{array}{c}\psi_1\\ \psi_2 \end{array}]:~\dR^{N+1}\to\dC^2$,
\be
\sL_\psi=\Tr\glh\dpsi^\da(\ri\pa_t\dpsi+\Delta\dpsi +V\dpsi)\grh\,,~~~\text{with}~~\osx\mapsto V(\osx)\in\dC^{2\ti2}\,,~V^\da=V\,,
\ee
leads to a $\dR$-valued Lagrangian functional $\cL$. Indeed
\[
\sL_\psi-\overline{\sL_\psi}=\ri\pa_t\Tr\glh\dpsi^\da\dpsi\grh+ \pa_{x\!_1}\Tr\glh\dpsi^\da(\pa_{x\!_1})\dpsi-(\pa_{x\!_1}\dpsi)^\da\dpsi\grh+\ldots+\pa_{x\!_N}\Tr\glh\dpsi^\da(\pa_{x\!_N})\dpsi-(\pa_{x\!_N}\dpsi)^\da\dpsi\grh\,.
\]
The $\sL_\psi$ of (2.27) leads to the Schr\"odinger equation for a particle with spin $\frac12$.
\end{examps}

\newpage

\section{ Free Gauge Fields}
The 'field variables' to be considered in this section are smooth functions
\be
\os{\aA}:\dR^N\,\to\,\underbrace{\dC^{c\ti c}\ti\cdots\ti\dC^{c\ti c}}_{N~ \text{times}}\,:\,\osx\mapsto\os{\aA}(\osx)=\col[\aA_1(\osx),\ldots,\aA_\mu(\osx),\ldots,\aA_N(\osx)]\,,
\ee
with $\aA_\mu(\osx)\,\in\,\gtg$, with $\gtg\subset\dC^{c\ti c}$ some fixed {\bf real} Lie algebra.
\footnote{In physics textbooks one often denotes $\ri\aA_\mu$, instead of $\aA_\mu$, cf. [DM]. For resemblance with Electromagnetism, I suppose. Because of $\gothic{u}(1)=
\ri\dR$ ?  To this author the factor $\ri$ is not convenient in all other cases.
}
This means that $\gtg$ is a $\dR$-linear subspace in $\dC^{c\ti c}$ which is not necessarily
a $\dC$-linear subspace.
On $\gtg$ we impose    the usual 'commutator'-Lie product
\[
\llh A_\mu\lk A_\nu\rlh=\big(A_\mu A_\nu-A_\nu A_\mu\big)\,.
\]

Important examples are   matrix Lie Algebras of type
\[
\gtg_J\,=\,\{\,X\in\dC^{r\ti r}\,\big|\,X^\dagger J+JX=0\,\}\,,~~~~~~\text{with fixed invertible}~~J\in\dC^{r\ti r}\,.
\]
Note that $\gtg_J$ is always a $\dR$-linear subspace in $\dC^{r\ti r}$, but not necessarily $\dC$-linear.

However: $\{J^{-1}=J^\dagger\}~\Rightarrow~\{X\in\gtg_J~\Rightarrow~X^\dagger\in\gtg_J\}$. 
\\[1mm]
Next, by $\proj: \dC^{c\ti c}\to\gtg$, we denote the {\bf real} orthogonal projection with respect to the {\bf real} inner product $X,Y\,\mapsto\,\Re\Tr[X^\dagger Y]$.

\begin{rems}\ \\
Consider $\dC^{c\ti c}$ as a {\bf real} vector space with standard {\bf real} inner product $X,Y\,\mapsto\,\Re\Tr[X^\dagger Y]$.

By $\proj: \dC^{c\ti c}\to\gtg$, we denote the {\bf real} orthogonal projection with respect this inner product.
\begin{itemize}
\item The Hermitean conjugation map $X\mapsto X^\dagger$ is $\dR$-linear symmetric and orthogonal.

\item If $\forall\,X\in\gtg: X^\da\in\gtg$, in short $\gtg^\da=\gtg$, it follows that $\forall X\in\dC^{c\ti c}: \proj (X^\dagger)=(\proj X)^\dagger$.

\item For fixed $K,L\in\dC^{c\ti c}$ the mapping $X\mapsto KX^\da L$ is $\dR$-linear. Its $\dR$-adjoint is $Y\mapsto LY^\da K$.

\item For any fixed invertble  $J\in\dC^{c\ti c}$ the mapping
\be
\sQ_J\,:\,\dC^{c\ti c}\to\dC^{c\ti c}~:~X\mapsto\sQ_J X \,=\,\frac12(X-J^{-1}X^\dagger J)\,,
\ee
is a $\dR$-linear mapping which reduces to the identity map when restricted to $\gtg_J$.

\item $\sQ_J$ is a $\dR$-linear projection on $\gtg_J$ iff $J=J^\dagger$.

\item $\sQ_J$ is a $\dR$-linear orthogonal projection on $\gtg_J$ if $J=J^{-1}=J^\dagger$.
\\ In this special case $\sQ_J=\proj$, with $\gtg=\gtg_J$.

\item If we modify the standard real inner product on $\dC^{c\ti c}$ to $X,Y\,\mapsto\,\Re\Tr[X^\dagger J^2 Y]$, the projection $\sQ_J$ is orthogonal
iff $J=J^\dagger$.
\end{itemize}
\end{rems}
{\bf Proof}

$\bullet$ $\Re\Tr[(X^\dagger)^\dagger Y]=\Re\Tr[XY]=\Re\Tr[X^\dagger (Y^\dagger)]$. Also $\Re\Tr[(X^\dagger)^\dagger(Y^\dagger)]=\Re\Tr[(X)^\dagger(Y)]$.

$\bullet$ Since $\gtg$ is supposed to be an invariant subspace for $X\mapsto X^\da$ and the latter is symmetric, also $\gtg^\perp$ is invariant.

$\bullet$ $\Re\Tr[(KX^\da L)^\da Y]=\Re\Tr[KX^\da L Y^\da]=\Re\Tr[X^\da(LY^\da K)]\,$.

$\bullet$ For $X\in\gtg$ holds $(I-\sQ_J) X=0\,$, iff $X\in\gtg\,$.

$\bullet$ $Q_J^2=Q_J$ iff $J=J^\da\,$. 



$\bullet\bullet$ $\frac12\Re\Tr[(X-J^{-1}X^\da J)^\da J^2 Y]=\frac12\Re\Tr[X^\da J^2 Y]-\frac12\Re\Tr[X^\da J^2(J^{-1}Y^\da J^{\da2}J^{-1})]\,$.

The 2nd term equals $-\frac12\Re\Tr[X^\da J^2(J^{-1}Y^\da J]\,$, for all $X,Y$, iff $J=J^\da\,$. \eindebewijs

\ \\[1mm]
Associated with  $\os{\aA}$, cf. (3.1), we introduce  covariant-type partial derivatives

$\nabla_\mu^A\,,1\leq\mu\leq N$ of functions $U\in \ms{C}^\infty(\dR^N\!:\!\dC^{c\ti c})$ by
\be
\nabla_\mu^A U\,=\,\pa_\mu U -\llh\cA_\mu\lk U\rlh\,=\,\pa_\mu U-\ad_{\cA_\mu}U\,.
\ee
One has the Leibniz-type rules
\be
\begin{array}{c}
\nabla_\mu^A (UV)\,=\,(\nabla_\mu^A U)V\,+\,U(\nabla_\mu^A V)\,,\hspace{1mm}\\[1.5mm]
\Tr\big[ U(\nabla_\mu^A V)\big]~=~\pa_\mu\Tr\big[ UV\big]-\Tr\big[(\nabla_\mu^A U)V\big]\,.
\end{array}
\ee
Note that if $U\in \ms{C}^\infty(\dR^N\!:\!\gtg)$ then also $\nabla_\mu^A U\in \ms{C}^\infty(\dR^N\!:\!\gtg)$.

\ \\[1mm]

Next, as in section 1, for given $\cA_\mu,\cA_\nu\in\fg\,,1\leq\mu,\nu\leq N$, define
\be
\cF_{\mu\nu}=\pa_\mu\cA_\nu-\pa_\nu\cA_\mu-\llh\cA_\mu\lk\cA_\nu\rlh~\in~\fg\,,
\ee
to which Theorem 1.2 applies.

For the construction of  a $\dR$-valued Lagrangian density $\sG_A$ for the Gauge field(s) $\os{\aA}$ we again employ a proto Lagrangian $\sG$, which is now an analytic function of $N(N-1)$ complex-matrix variables 
 and just smooth in $N$ real variables:
\be
\sG~:~\underbrace{\dC^{c\ti c}\ti\cdots\ti\dC^{c\ti c}}_{\frac12 N(N-1)~ \text{times}}~\ti~\underbrace{\dC^{c\ti c}\ti\cdots\ti\dC^{c\ti c}}_{\frac12 N(N-1)~ \text{times}}~\ti~\dR^N~~\to~~\dC\,.
\ee
The  1st set of entries to this function is labeled by the ordered pairs $(\mu\nu)\,,\,1\leq\mu<\nu\leq N\,$. The 2nd set of entries is labelled by the ordered triple $(\th\rho\star)\,,\,1\leq\th<\rho\leq N\,$. We denote
\[
\{~\ldots ,P_{\mu\nu},\ldots;\ldots, Q_{\th\rho\star},\ldots;\osx\}~\mapsto~\sG(~\ldots P_{\mu\nu},\ldots;\ldots Q_{\th\rho\star},\ldots;\osx)\in\dC\,,
\]
with $1\leq\mu<\nu\leq N$ and $1\leq\th<\rho\leq N$. 
The 3 bunches  of variables get their corresponding partial derivatives denoted by, respectively, cf. (2.4),
\[
\sG^{(\dg{\mu\nu})}(\ldots,P_{\th\rho},\ldots;\ldots,Q_{\th\rho\star},\ldots;\osx)\,,~~\sG^{(\dg{\th\rho\star})}(\ldots,P_{\th\rho},\ldots;\ldots,Q_{\th\rho\star},\ldots;\osx)\,,~~\sG^{(\dg{\nabla})}\,.
\]
Let the Lie algebra $\gtg$ be fixed. On $\sG$ we put the  condition, take $Q_{\th\rho\star}=P_{\th\rho}^\dagger$,
\be~~\forall~\{P_{\mu\nu}\}_{1\leq\mu<\nu\leq N}\subset\gtg\,~~\forall\osx\in\dR^N~:~\sG(\ldots,P_{\mu\nu},\ldots;\ldots,P_{\th\rho}^\dagger,\ldots;\osx)\,\in\,\dR\,.
\ee
 The   Lagrangian density we want to consider is found by replacing $P_{\mu\nu}\to\aF_{\mu\nu}\,,Q_{\th\rho\star}\to\aF_{\th\rho}^\dagger$,
\be
\osx~~\mapsto~~\sG_A(\osx)=\sG(~\ldots, \aF_{\mu\nu}(\osx),\ldots;\ldots, \aF^\dagger_{\th\rho}(\osx),\ldots;\osx\,)\,\in\dR.
\ee
Note that if $\gtg=\gtg_J$, for some fixed $J\in\dC^{c\ti c}$, we have $\aF_{\th\rho}^\dagger=-J\aF_{\th\rho}J^{-1}\,,\th<\rho$.

As in the previous section, a corresponding useful notation is
\be
\osx~~\mapsto~~\sG^{(\dg{\mu\nu})}_A(\osx)=\sG^{(\dg{\mu\nu})}(~\ldots, \aF_{\mu\nu}(\osx),\ldots;\ldots, \aF^\dagger_{\th\rho}(\osx),\ldots;\osx\,)\,\in\dC^{c\ti c}.
\ee
The Lagrangian density $\sG_A$ depends on the field variables $\osx\mapsto\cA_\mu(\osx)\,,1\leq\mu\leq N$, and their derivatives. All being functions in a vectorspace over $\dR$.
In the important special case $\gtg=\gtg_J$ the hermitean conjugate notation of the field variables $\aA_\mu$ need not even occur.

Finally, note that, because of (2.11) and (3.8), we have
\be
\sG^{(\dg{\th\rho\star})}_A(\osx)=(\sG^{(\dg{\th\rho})}_A)^\dagger(\osx)\,,~~1\leq\th<\rho\leq N\,.
\ee

\begin{notation}
 In order to visually simplify the formulae to come, it
 is useful to extend the set of functions $\sG^{(\dg{\mu\nu})}_A$, cf.(3.9), to  'full' labels $1\leq\mu,\nu\leq N$ in the following way,
\be
~\hat{\sG}^{(\dg{\mu\nu})}_A=\left\{\begin{array}{cl}
\sG^{(\dg{\mu\nu})}_A & \text{if}~~ 0\leq\mu<\nu\leq N\,, ~\text{ as before}, \\[1mm]
0                     & \text{if}~~ \mu=\nu\,,                                 \\[1mm]
-\sG^{(\dg{\nu\mu})}_A & \text{if}~~ 0\leq\nu<\mu\leq N\,.\\
\end{array}\right.
\ee
\end{notation}

\bt
Fix a matrix Lie algebra $\gtg\subset\dC^{c\ti c}$. Consider the Lagrangian density $\sG_A$ of (3.8).\\[1mm]
{\bf A.} The  Euler-Lagrange equations for the free gauge fields $\aA_{\mu}\,,1\leq\mu\leq N$, with values in the Lie algebra $\gtg\subset\dC^{c\ti c}$, read
\be
\sum_{\mu=1}^N\proj\Big(\Big(\nabla^A_\mu\big([\proj\hat{\sG}^{(\dg{\mu\ka\star})}_A]^\da\big)\Big)^\da\Big)~=~0\,,~~~1\leq\ka\leq N\,,
\ee
with $\nabla^A_\mu$ as in (3.3).\\[1mm]
{\bf B.} In the special case $\gtg^\da=\gtg$ the Euler-Lagrange equations simplify to
\be
\sum_{\mu=1}^N\Big(\nabla_\mu^A\proj\hat{\sG}^{(\dg{\mu\ka})}_A\Big)~=~0\,,~~~1\leq\ka\leq N\,.
\ee
{\bf C.} If we take $\gtg=\gtg_J$, with $J=J^\dagger=J^{-1}$, the latter becomes
\be
\sum_{\mu=1}^N\nabla_\mu^A\Big(\sQ_J[\hat{\sG}^{(\dg{\mu\ka})}_A]\Big)\,=\,0\,,~~~~1\leq\ka\leq N\,,
\ee
where $\sQ_J \mZ=\frac12\mZ-\frac12 J\mZ^\da J\,,~\mZ\in\dC^{c\ti c}$.
\et

{\bf Proof} \\
{\bf A.} In order to calculate the (directional) derivatives of the Lagrangian functional $\cG=\int\sG_A\,\rd\osx$ with respect to the free gauge fields $\aA_\ka\,,1\leq\ka\leq N\,$, we first expand a perturbation of $\osx\mapsto\aF_{\mu\nu}(\osx)$ by substitution of the gauge fields $\osx\mapsto \aA_\mu(\osx)+\vep\del_{\mu\ka}\aH(\osx)\,,\vep\in\dR$,
\[
\aF_{\mu\nu;\vep,\ka}=\glh\pa_\mu(\cA_\nu+\vep\del_{\nu\ka}\cH)-\pa_\nu(\cA_\mu+\vep\del_{\mu\ka}\cH)
-\llh\cA_\mu+\vep\del_{\mu\ka}\cH \lk \cA_\nu
+\vep\del_{\nu\ka}\cH\rlh\,\grh = \hspace{1.4cm}
\]
\[
=\,\glh\pa_\mu\cA_\nu-\pa_\nu\cA_\mu-\llh\cA_\mu \lk \cA_\nu\rlh\grh\,+\,
\vep\,\del_{\nu\ka}\glh\pa_\mu\cH-\llh\cA_\mu \lk \cH\rlh\,\grh\,
-\vep\,\del_{\mu\ka}\glh\pa_\nu\cH-\llh\cA_\nu \lk \cH\rlh\,\grh
\,=
\]
\[
\hspace{3cm}=\,\aF_{\mu\nu}+\vep\del_{\nu\ka}\nabla^A_\mu \aH-\vep\del_{\mu\ka}\nabla^A_\nu \aH\,.
\]

Consider the expansion
\[
\sG(\ldots,\aF_{\mu\nu;\vep,\ka}\,,\ldots;\ldots,\aF_{\th\rho;\vep,\ka}^\dagger,\ldots;\osx)\,-\,\sG(\ldots,\aF_{\mu\nu}\,,\ldots;\ldots,\aF_{\th\rho}^\dagger,\ldots;\osx)\,=\hspace{1cm}
\]
\[
  =\vep\sum_{1\leq\mu<\nu\leq N}\Tr\glh[\sG^{(\dg{\mu\nu})}_A][\del_{\nu\ka}\nabla^A_\mu \aH-\del_{\mu\ka}\nabla^A_\nu \aH\,]\,+ \hspace{7cm}
\]
\[
 \hspace{4cm} +\, \vep\sum_{1\leq\th<\rho\leq N} \Tr\glh[\sG^{(\dg{\th\rho\star})}_A][\del_{\rho\ka}\nabla^A_\th \aH-\del_{\th\ka}\nabla^A_\rho \aH_\ka\,]^\dagger\grh~+\mathscr{O}(\vep^2)\,=
\]
\[
  =\frac\vep2\sum_{\mu,\,\nu=1}^N\Tr\glh[\hat{\sG}^{(\dg{\mu\nu})}_A][\del_{\nu\ka}\nabla^A_\mu \aH-\del_{\mu\ka}\nabla^A_\nu \aH\,]\,+ \hspace{7cm}
\]
\[
 \hspace{4cm} +\, \frac\vep2\sum_{\th,\,\rho=1}^N \Tr\glh[\hat{\sG}^{(\dg{\th\rho\star})}_A][\del_{\rho\ka}\nabla^A_\th \aH-\del_{\th\ka}\nabla^A_\rho \aH\,]^\dagger\grh~+\mathscr{O}(\vep^2)\,=
\]

\[
  =\frac\vep2\sum_{\mu=1}^N\Tr\glh[\hat{\sG}^{(\dg{\mu\ka})}_A][\nabla^A_\mu \aH\grh-\frac\vep2\sum_{\nu=1}^N\Tr\glh[\hat{\sG}^{(\dg{\ka\nu})}_A][\nabla^A_\nu \aH\,]\grh\,+ \hspace{5cm}
\]
\[
 \hspace{4cm} +\, \frac\vep2\sum_{\th=1}^N \Tr\glh[\hat{\sG}^{(\dg{\th\ka\star})}_A][\nabla^A_\th \aH\grh\,-
 \frac\vep2\sum_{\rho=1}^N \Tr\glh[\hat{\sG}^{(\dg{\ka\rho\star})}_A][\nabla^A_\rho \aH\,]^\dagger\grh~+\mathscr{O}(\vep^2)\,=
\]
\[
  =\vep\sum_{\mu=1}^N\Tr\glh[\hat{\sG}^{(\dg{\mu\ka})}_A][\nabla^A_\mu \aH\grh~+~\vep\sum_{\mu=1}^N \Tr\glh[\hat{\sG}^{(\dg{\mu\ka\star})}_A][\nabla^A_\mu \aH\,]^\dagger\grh~+\mathscr{O}(\vep^2)\,= \hspace{3cm}
\]
\[
  =2\vep\Re\sum_{\mu=1}^N\Tr\glh[\hat{\sG}^{(\dg{\mu\ka\star})}_A]^\da [\nabla^A_\mu \aH]\grh~+~\mathscr{O}(\vep^2)\,=\, 2\vep\Re\sum_{\mu=1}^N\Tr\glh[\proj\hat{\sG}^{(\dg{\mu\ka\star})}_A]^\da [\nabla^A_\mu \aH]\grh~+~\mathscr{O}(\vep^2)\,=
\]
\[
=\, -2\vep\Re\sum_{\mu=1}^N\Tr\glh\nabla^A_\mu\big([\proj\hat{\sG}^{(\dg{\mu\ka\star})}_A]^\da\big) \aH\grh~+~\sum_{\mu=1}^N\pa_\mu(\ldots)
\,+\mathscr{O}(\vep^2)\,=
\]

\be
=\, -2\vep\Re\sum_{\mu=1}^N\Tr\glh\,\Big(\proj\Big(\Big(\nabla^A_\mu\big([\proj\hat{\sG}^{(\dg{\mu\ka\star})}_A]^\da\big)\Big)^\da\Big)\Big)^\da \aH\grh~+~\sum_{\mu=1}^N\pa_\mu(\ldots)
\,+\mathscr{O}(\vep^2)\,.
\ee
In this derivation we used, respectively, the antisymmetry $\mu\leftrightarrow\nu$ of $[\hat{\sG}^{(\dg{\mu\nu})}_A]$ and $[\del_{\nu\ka}\nabla^A_\mu \aH-\del_{\mu\ka}\nabla^A_\nu \aH\,]$, the Leibniz rule(3.4),  the fact that $\Re\Tr\glh\big(\ldots\big)^\da\aH\grh$ expresses the real inner product on $\dC^{c\ti c}$ and $\proj$ the real orthogonal projection on $\gtg$.

Also properties like $\Tr[AB]=\Tr[BA]\,$, $\Tr[A\llh B\lk C\rlh]=\Tr[\llh A\lk B\rlh C]\,$ play a crucial role.

The result now follows by the usual variational practices.

{\bf B.} If $\gtg^\da=\gtg$ the real linear mappings $\{.\}^\da$ and $\proj$ commute, which greatly simplifies the result of A.

{\bf C.} Use Remarks 3.1.
\eindebewijs
\ \\[1mm]

\begin{examp}\ \\ 
{\bf A.} For convenience we restrict to Lie-algebras with property $\gtg^\da=\gtg$. We will consider general Lagrangians which are (real) quadratic in $\aF_{\mu\nu}$. Here, in our summation expressions, we
write $\mu<\nu$ instead of $1\leq\mu<\nu\leq N$. Start from the proto Lagrangian
\be
\sG=\sum_{\mu<\nu\,,\,\th<\rho} h_{(\mu\nu)(\th\rho)}\Tr[P_{\mu\nu}Q_{\th\rho\star}]\,~~\text{with}~~\bs{h_{(\mu\nu)(\th\rho)}}=h_{(\th\rho)(\mu\nu)}\in\dC\,.
\ee
Note
\[
\sum_{\mu<\nu,\th<\rho} h_{(\mu\nu)(\th\rho)}\Tr[P_{\mu\nu}P_{\th\rho}^\dagger]\,\in\,\dR\,.
\]
For the derivatives of $\sG$ we find, 
\[
~~~\sG^{(\dg{\mu\nu})}(\ldots,P_{\mu\nu},\ldots;\ldots,Q_{\th\rho\star},\ldots\,)=\sum_{\al<\bet} h_{(\mu\nu)(\al\bet)}Q_{\al\bet\star}
\]
\[
~~~\sG^{(\dg{\th\rho\star})}(\ldots,P_{\mu\nu},\ldots;\ldots,Q_{\th\rho\star},\ldots\,)=\sum_{\al<\bet} h_{(\al\bet)(\th\rho)}P_{\al\bet}
\]
If we take $Q_{\th\rho\star}=P_{\th\rho}^\dagger\,$, one easily checks (3.8),
\[
{\sG^{(\dg{\mu\nu})}}^\dagger(\ldots,P_{\mu\nu},\ldots;\ldots,P_{\th\rho}^\da,\ldots\,)=\sum_{\al<\bet} \bs{h_{(\mu\nu)(\al\bet)}}P_{\al\bet}=\sum_{\al<\bet} h_{(\al\bet)(\mu\nu)}P_{\al\bet}=
\sG^{(\dg{\mu\nu\star})}\,.
\]
The Lagrangian density
\be
\sG_A=\sum_{\mu<\nu,\,\th<\rho} h_{(\mu\nu)(\th\rho)}\Tr[\cF_{\mu\nu}\cF_{\th\rho}^\dagger]\,,
\ee
can now be put in (3.13) to find the  Euler-Lagrange equations. Note however, that $\proj$ cannot be put 'through' the $h_{(\mu\nu)(\th\rho)}$ if those are non-real numbers!

So, let us  restrict to $\gtg^\da=\gtg$ {\`a}nd $h_{(\mu\nu)(\th\rho)}\in\dR$. Anti-symmetrize  $h_{(\mu\nu)(\th\rho)}$ to full labels:
\[
~\hat{h}_{(\mu\nu)(\th\rho)}=\left\{\begin{array}{cl}
h_{(\mu\nu)(\th\rho)} & \text{if}~~ \mu<\nu\,,\th<\rho ~~\text{or} ~~\mu>\nu\,,\th>\rho\\[1mm]
0                     & \text{if}~~ \mu=\nu\,~\text{and/or} ~\th=\rho                                \\[1mm]
-h_{(\nu\mu)(\th\rho)} & \text{if}~~ \mu>\nu\,,\th<\rho\,\\[1mm]
-h_{(\mu\nu)(\rho\th)} & \text{if}~~ \mu<\nu\,,\th>\rho\,\\
\end{array}\right.
\]
In this special case
\[
\hat{\sG}^{(\dg{\mu\nu})}_A~=~\frac12\sum_{\al,\bet=1}^N \hat{h}_{(\mu\nu)(\al\bet)}\aF_{\al\bet}^\da~,
\]
 and, since $\aF_{\al\bet}^\da\in\gtg$, the E-L-equations (3.13) become
\be
\frac12\sum_{\al,\bet=1}^N\sum_{\mu=1}^N \hat{h}_{(\mu\ka)(\al\bet)}\Big(\pa_\mu\aF_{\al\bet}^\da-\llh\aA_\mu\lk\aF_{\al\bet}^\da\rlh\Big)~=~0~,~~~1\leq\ka\leq N\,.
\ee

\ \\[1mm]

{\bf B.} For gauge fields on Minkowski space, with coordinates $x^0,x^1,x^2,x^3$ and

metric $[g^{\mu\nu}]=\mr{diag}(1,-1,-1,-1)$, one usually takes, cf. [DM],
\[
h_{(\mu\nu)(\al\bet)}=g^{\mu\al}g^{\nu\bet}=(-1)^{1+\del_{\mu0}}\del_{\mu\al}(-1)^{1+\del_{\nu0}}\del_{\nu\bet}=(-1)^{\del_{\mu0}+\del_{\nu0}}\del_{\mu\al}\del_{\nu\bet}\,.
\]
Hence
\[
\hat{h}_{(\mu\ka)(\al\bet)}~=~\sgn(\ka-\mu)\,\sgn(\bet-\al)\,(-1)^{\del_{\mu0}+\del_{\ka0}}\del_{\mu\al}\del_{\ka\bet}\,.
\]
In this special case the Lagrangian density (3.17) reads
\be
\sG_A=\sum_{0\leq\mu<\nu\leq 3}\,(-1)^{\del_{\mu0}+\del_{\nu0}}\Tr\big[\aF_{\mu\nu}  \aF_{\mu\nu}^\dagger \,\big]\,.
\ee
The corresponding Euler-Lagrange equations are 
\be
\sum_{\mu=0}^3 (-1)^{\del_{\mu0}+\del_{\ka0}}\nabla_\mu^A\,\aF_{\mu\ka}^\da=0\,,~~~0\leq\ka\leq3\,.
\ee

For $\dim\gtg=1$ the term $\ad_{\aA_\mu} \aF_{\mu\ka}^\dagger$ vanishes. This simplification, viz. $\nabla_\mu^A=\pa_\mu~$, leads to  standard electromagnetism in Minkowski space.
Indeed, if we put $\aA_0^\dagger=-\Phi$ and $\col[\aA_1^\dagger\,,\aA_2^\dagger\,,\aA_3^\dagger]=\os{A}$, then (3.20) turns into Maxwell's equations 'in potential form'
\be
\left\{\begin{array}{rcc}
\paf{}{t}\div\os{A} +\Delta \Phi  & =& 0\\[2mm]
\paft{}{t}\os{A}-\Delta\os{A}+\grad\big(\paf{}{t}\Phi+\div\os{A}\big) &=&\os{0}
\end{array}
\right.
\ee
If the pair $\os{A},\os{B}$ satisfies (3.21), then the pair $\os{E}=-\paf{\os{A}}{t}-\grad\Phi\,,\,\os{B}=\rot\os{A}\,$, satisfies the classical Maxwell equations.

Finally, imposing the 'Lorenz-Gauge' $~\paf{}{t}\Phi+\div\os{A}=0$, we find the usual wave equations $\pa_t^2\Phi-\Delta\Phi=0\,,\,\pa_t^2\os{A}-\Delta\os{A}=\os{0}\,$.
For more details see Appendix B.
\end{examp}

\section{Noether  Fluxes}
'Infinitesimal symmetries' of the Lagrangian density $\sL$ lead to local conservation laws for the solutions of the Euler Lagrange equations. So we are told by Emmy Noether's famous theorem. First we have a short look at the needed concepts as formulated within our special (simple) context.
\begin{defn}
A {\em Conservation Law} or {\em Noether Flux} is a vectorfield on $\dR^N$, with components $\sV^\mu_\psi\,,\,1\leq\mu\leq N$, which arise from a set of functions of Proto-Lagrangian type, $\sV^\mu\,,\,1\leq\mu\leq N\,$\,, cf. (2.13),
such that for all solutions $\dpsi$ of the Euler Lagrangian system, cf. Th 2.4, we have
\be
\sum_{\mu=1}^ N \paf{}{x^\mu}\sV_\psi^\mu(\osx)=0\,,~~~~~~\text{where}~~\sV_\psi^\mu(\osx)=\sV^\mu(\dpsi(\osx),\dpsi^\dagger(\osx),\dpsi_{,\mu}(\osx),\dpsi^\dagger_{,\mu}(\osx),\osx)\,.
\ee
\end{defn}
A conservation law can be  named  'trivial' for several reasons: It may happen that for {\em all} solutions $\dpsi$ the fluxes $\sV^\mu_\psi=0$. Another reason for triviality occurs if for all functions $\dpsi$, whether they are solutions or not,
the identity (4.1) is satisfied. For example if the components $\sV^\mu_\psi$ arise from the curl of an arbitrary vector field depending on $\dpsi$.\\[1mm]

Two types of symmetries will be considered here: '{Internal symmetries}'
and '{External symmetries}'. They can be formulated in terms  of the {\em proto-Lagrangian} only.\\
External symmetries regard transformations 
of the spatial variables $\osx$. We restrict to {\em affine transforms}. 

\begin{defn} ({\bf Internal symmetries})\\
A set of linear mappings $\rK\,,\rL_\mu^\la\,:\dC^{r\ti c}\to\dC^{r\ti c}\,,1\leq\la,\mu\leq N$, is said to generate an internal (local) symmetry of the proto-Lagrangian $\sL$ if for all $\mP,\mQ_\mu\in\dC^{r\ti c}$, all $\osx\in\dR^N$,
and  $s\in\dR$, $|s|$ small, one has                                          
\[ \sL(e^{s\rK}\mP;(e^{s\rK}\mP)^\dagger;\ldots e^{s\rL^\la_\mu}\mQ_\la\ldots;\ldots(e^{s\rL^\la_\mu}\mQ_\la)^\dagger\ldots;\,\osx) =\hspace{3cm}
\]
\be
~\hspace{5cm}=\sL(\mP;\mP^\dagger;\ldots \mQ_\mu\ldots;\ldots\mQ_\mu^\dagger\ldots;\,\osx)+\ms{O}(s^2)\,,
\ee
\end{defn}
In many cases the $\rK\,,\rL^\la_\mu$ are realized by left and/or right multiplication with some fixed matrices in $\dC^{r\ti r}$ or $\dC^{c\ti c}$.\\[0.5mm]
Many times there is a special type of internal symmetry which is related to a linear mapping $A:\dR^N\to\dR^N$ in the 'outside world',
\[
 \sL(\mP;\mP^\dagger;\ldots (e^{sA})^\la_\mu\mQ_\la\ldots;\ldots((e^{sA})^\la_\mu\mQ_\la)^\dagger\ldots;\,\osx) =\hspace{3cm}
\]
\be
~\hspace{5cm}=\sL(\mP;\mP^\dagger;\ldots \mQ_\mu\ldots;\ldots\mQ_\mu^\dagger\ldots;\,\osx)+\ms{O}(s^2)\,,
\ee
\begin{defn} ({\bf External symmetries})\\
The affine mapping $\osx\mapsto -s\osa+e^{sA}\osx$ on $\dR^N$, where $\osa\in\dR^N$ and $A:\dR^N\to\dR^N$, a linear mapping, is said to generate an external (local) symmetry of the proto-Lagrangian $\sL$ if for all $\mP,\mQ_\mu\in\dC^{r\ti c}$, all $\osx\in\dR^N$,
and  $s\in\dR$, $|s|$ small, one has
\[
 \sL(\mP;\mP^\dagger;\ldots \mQ_\mu\ldots;\ldots(\mQ_\mu)^\dagger\ldots;\,-s\osa+e^{sA}\osx) =\hspace{3cm}
\]
\be
~\hspace{5cm}=\sL(\mP;\mP^\dagger;\ldots \mQ_\mu\ldots;\ldots\mQ_\mu^\dagger\ldots;\,\osx)+\ms{O}(s^2)\,.
\ee
\end{defn}
\begin{remarks} \ \vspace{-8mm}
\begin{itemize}
\item The order constant in $\ms{O}(s^2)$ may depend on all independent variables of $\sL$.

\item If in (4.2)-(4.4) exponents like $e^{s\rK}$ are replaced by $\rI +s\rK$ we get equivalent conditions. However in many practical applications the terms $\ms{O}(s^2)$
are identically zero if exponentials are used.
\item Local symmetry (4.4) implies
\[
\sL^{(\dg{\nabla})}(\mP;\mP^\dagger;\ldots \mQ_\mu\ldots;\ldots\mQ_\mu^\dagger\ldots;\,\osx)\cdot(A\osx-\osa)=0\,.
\]
\end{itemize}
\end{remarks}

We now first consider two types of  conservation laws in connection with affine transformations  in space.

 For any vector $\os{a}\in\dR^N$ we define the {\em Translation operator} $\mathbf{T}_{\os{a}}$ by
\[
\mathbf{T}_{\os{a}}\dpsi(\osx)=\dpsi(\osx-\os{a}).
\]
For any matrix $A\in\dR^{N\times N}$ we define the {\em dilation operator} $\mathbf{R}_A$ by
\[
\mathbf{R}_A\dpsi(\osx)=\dpsi(e^{A}\osx).
\]
\bt
Suppose that, for some $\rK:\dC^{r\ti c}\to\dC^{r\ti c}$ and some $\osa\in\dR^N$, the proto-Lagrangian $\cL$ has internal local symmetry (4.2) with $\rL^\la_\mu=\del^\la_\mu\rK$ and
external local symmetry (4.4) with $A=O$. Then for any solution $\dpsi$ of the Euler-Lagrange system one has the conservation law
\be
\sum_{\mu=1}^N\frac{\pa}{\pa x^\mu}\Big\{\Tr\glh[\sL^{\dg{(\mu)}}_\psi]\cdot(\rK\dpsi- a^\la\pa_\la\dpsi)+ [\sL^{\dg{(\mu\star)}}_\psi]\cdot(\rK\dpsi- a^\la\pa_\la\dpsi)^\dagger\grh+a^\mu\sL_\psi\Big\}=0~.~~~~~~
\ee
\et
{\bf Proof:}
By $\cong$ we mean equality up to a term   $\ms{O}(s^2)$. We study 
\[
\sL\big(\, e^{s\rK}\mathbf{T}_{s\os{a}}\dpsi ,\,\mathbf{T}_{s\os{a}}\dpsi^{\dagger} e^{s\rK^{\dagger}},\,\pa_\mu[e^{sK}\mathbf{T}_{s\os{a}}\dpsi ],\,\pa_\mu[\mathbf{T}_{s\os{a}}\dpsi^{\dagger} e^{sK^{\dagger}}],\,\osx-s\osa\big)\,.
\]
With our conditions it can be written
\[
\sL(e^{s\rK}\dpsi(\os{x}-s\os{a});( e^{s\rK}\dpsi(\os{x}-s\os{a}))^\dagger;\ldots\pa_\mu e^{s\rK}\dpsi(\os{x}-s\os{a})\ldots;\ldots\pa_\mu( e^{s\rK}\dpsi(\os{x}-s\os{a}))^\dagger\ldots;\,\osx-s\osa) \cong
\]
\[
\cong\sL(\dpsi(\os{x}-s\os{a});\dpsi(\os{x}-s\os{a})^\dagger;\ldots \dpsi_{,\mu}(\os{x}-s\os{a})\ldots; \ldots\dpsi_{,\mu}(\os{x}-s\os{a})^\dagger\ldots;\,\osx-s\osa) =
\]
\be
=\sL_\psi(\os{x}-s\os{a})=(\mathbf{T}_{s\os{a}}\sL_\psi)(\osx)\,.
\ee
Differentiate the first line of this at $s=0$ and use $\sL^{(\dg{\nabla})}\!\cdot\!\osa=0\,$,
\[
\Tr\big\{[\sL^{\dg{(o)}}_\psi](\rK\dpsi- a^{\la}\pa_\la\dpsi)+ [\sL^{\dg{(o\star)}}_\psi](\dpsi^\dagger \rK^\dagger -a^{\la}\pa_\la\dpsi^\dagger)+\hspace{5cm}
\]
\be \hspace{2cm}
+[\sL^{\dg{(\mu)}}_\psi](\rK\pa_\mu\dpsi- a^\la\pa_\la\pa_\mu\dpsi)+ [\sL^{\dg{(\mu\star)}}_\psi](\pa_\mu\dpsi^\dagger \rK^\dagger -a^\la\pa_\la\pa_\mu\dpsi^\dagger)
\big\}.
\ee
If $\dpsi$ is a solution we use (2.16) and replace $[\sL^{\dg{(o)}}_\psi]$ by $\frac{\pa}{\pa x^\mu}[\sL^{\dg{(\mu)}}_\psi]$, etc.
Now (4.7) can be written as a divergence, which constitutes  the left hand side of (4.5), apart from the last term inside $\{~~~\}$. Together with the derivative $a^\la\pa_\la\sL_\psi=\pa_\mu(a^\mu\sL_\psi)$
at $s=0$ of the final line of (4.6) we arrive at the  wanted conserved current (4.5).
\hfill{$\blacksquare$}\\[3mm]

\begin{examp}
 Let $\Gamma^\mu$ and $M$ be constant complex matrices with $\Gamma^{\mu\dagger}=\Gamma^\mu$ and $M=-M^\dagger$. Then the Lagrangian density
\be
\sL_\psi=\Tr\big\{\ri\dpsi^\dagger \Gamma^\mu\pa_\mu\dpsi+\dpsi^\dagger M\dpsi \big\},
\ee
for $\dpsi:~\dR^N\to\dC^{r\times c}$ satisfies the condition of Theorem 4.1 for $K=O$ and all $\os{a}\in\dR^N$.
The conservation law reads
\be
\frac{\pa}{\pa x^\mu}\Trc{-a^{\lambda}\dpsi^\dagger\Gamma^\mu\pa_\la\dpsi+ a^{\mu}\dpsi^\dagger\Gamma^\lambda\pa_\la\dpsi
+a^\mu\dpsi^\dagger M\dpsi}\,=\,\frac{\pa}{\pa x^\mu}\Trc{-a^{\lambda}\dpsi^\dagger\Gamma^\mu\pa_\la\dpsi}\,=\,0.
\ee
This can be checked directly for solutions of the PDE: $~\Gamma^\mu\pa_\mu\dpsi+M\dpsi=0$. Observe that in this special case $\sL_\psi=0$ for solutions.\\
Also the Lagrangian of Example (2.5b), with {\em constant} matrices $K,\, M,\,\Gamma^\mu,\,\cA_\mu$ leads to conservation laws of this type.
\end{examp}

\bt
Suppose that, for some $\rK:\dC^{r\ti c}\to\dC^{r\ti c}$ and some $A\in\dR^{N\ti N}$ with $\Tr A=0$, the proto-Lagrangian $\cL$ has internal local symmetry (4.2)
with $\rL^\la_\mu=\rK+[A]_\mu^\la\rI$ and
external local symmetry (4.4) with $\osa=\os{0}$. Then for any solution $\dpsi$ of the Euler-Lagrange system one has the conservation law
\[
\sum_{\mu=1}^N\frac{\pa}{\pa x^\mu}\Big\{\Tr\glh[\sL^{\dg{(\mu)}}_\psi](\rK\dpsi(\osx)+A^\alpha_\beta\,x^\beta\,\dpsi_{,\alpha}(\osx))\,+\hspace{5cm}
\]
\be
\hspace{5cm}+\,[\sL^{\dg{(\mu\star)}}_\psi](\rK\dpsi(\osx)+A^\alpha_\beta\,x^\beta\,\dpsi_{,\alpha}(\osx))^\dagger\grh\,-\,A^\mu_\beta x^\beta\,\sL_\psi\Big\}\,=\,0\,.
\ee
\et
{\bf Proof:} We study
\[
\sL\big(\, e^{s\rK}\mathbf{R}_{sA}\dpsi\, ;\,\mathbf{R}_{sA}\dpsi^{\dagger}e^{s\rK^\dagger}\,;\ldots\pa_\mu[e^{s\rK}\mathbf{R}_{sA}\dpsi]\ldots\,;
\ldots\pa_\mu[\mathbf{R}_{sA}\dpsi^{\dagger}e^{s\rK^\dagger}]\,\ldots;e^{sA}\osx\,\big)\,.
\]
With our conditions it can be written,
\[
\sL(\dpsi(e^{sA}\osx);\dpsi(e^{sA}\osx)^\dagger;\ldots\pa_\mu \dpsi(e^{sA}\osx)\ldots;\ldots\pa_\mu\dpsi(e^{sA}\osx)^\dagger\ldots;\,e^{sA}\osx) \cong
\]
\[
\cong\sL(\dpsi(e^{sA}\osx);\dpsi(e^{sA}\osx)^\dagger;\ldots(e^{sA})^\la_\mu \dpsi_{,\la}(e^{sA}\osx)\ldots;\ldots(e^{sA})^\la_\mu\dpsi_{,\la}(e^{sA}\osx)^\dagger\ldots;\,e^{sA}\osx) \cong
\]
\[
\cong\sL(\dpsi(e^{sA}\osx);\dpsi(e^{sA}\osx)^\dagger;\ldots \dpsi_{,\mu}(e^{sA}\osx)\ldots;\ldots\dpsi_{,\mu}(e^{sA}\osx)^\dagger\ldots;\,e^{sA}\osx) \cong
\]
\[
\cong\sL(\dpsi(e^{sA}\osx);\dpsi(e^{sA}\osx)^\dagger;\ldots \dpsi_{,\mu}(e^{sA}\osx)\ldots;\ldots\dpsi_{,\mu}(e^{sA}\osx)^\dagger\ldots;\,e^{sA}\osx)\,=
\]
\be =\sL_\psi(e^{sA}\osx)\,=\,(\mathbf{R}_{sA}\sL_\psi)(\osx).
\ee
Differentiate the first line of this at $s=0$\, and use $\sL^{(\dg{\nabla})}\!\cdot\!A\osx=0\,$:
\[
\Tr\big\{[\sL^{\dg{(o)}}_\psi](\rK\dpsi(\osx)+A^\alpha_\beta\,x^\beta\,\dpsi_{,\alpha}(\osx))+ [\sL^{\dg{(\mu)}}_\psi]\pa_\mu(\rK\dpsi(\osx)+A^\alpha_\beta\,x^\beta\,\dpsi_{,\alpha}(\osx))\,+\hspace{5cm}
\]
\be \hspace{7mm}
+\,[\sL^{\dg{(o\star)}}_\psi](\rK\dpsi(\osx)+A^\alpha_\beta\,x^\beta\,\dpsi_{,\alpha}(\osx))^\dagger+ [\sL^{\dg{(\mu\star)}}_\psi]\pa_\mu(\rK\dpsi(\osx)+A^\alpha_\beta\,x^\beta\,\dpsi_{,\alpha}(\osx))^\dagger
\big\}.
\ee
If $\dpsi$ is a solution we use (2.16) and replace $[\sL^{\dg{(o)}}_\psi]$ by $\frac{\pa}{\pa x^\mu}[\sL^{\dg{(\mu)}}_\psi]$, etc.
Now (4.12) can be written as a divergence, which constitutes  the left hand side of (4.10), apart from the last term between $\{~~~\}$. Together with
the derivative at $s=0$ of the final line in (4.11):
$A^\mu_\beta \pa_\mu\sL_\psi=\pa_\mu(A^\mu_\beta x^\beta\,\sL_\psi)$, use $\Tr A=0$,
 we arrive at the  conserved current (4.10). \eindebewijs

\  \\[2mm]
Next we  deal with {\em\bf internal symmetries} only. They play a crucial role in Gauge theories. A simple case first.


\bt
Suppose that, for some linear $\rK:\dC^{r\ti c}\to\dC^{r\ti c}$  the proto-Lagrangian $\cL$ satisfies (4.2) with $\rL^\la_\mu=\del^\la_\mu\rK$. Then for any solution $\dpsi$ of the Euler-Lagrange system one has the conservation law


\be
\sum_{\mu=1}^N\frac{\pa}{\pa x^\mu}\Tr\big\{[\sL^{\dg{(\mu)}}_\psi]\rK\dpsi  + [\sL^{\dg{(\mu\star)}}_\psi] (\rK \dpsi)^\dagger\big\}=0~,~~~~~~
\ee
\et
{\bf Proof:} Calculate the derivative
\[
\dfrac{\pa}{\pa s}\sL\big(e^{s\rK}\dpsi ,\,(e^{s\rK}\dpsi)^\dagger,\,\pa_\mu[e^{s\rK}\dpsi],\pa_\mu [e^{s\rK}\dpsi]^\dagger,\,\osx\,\big),~~~~\mbox{at}~~~~s=0\,.
\]
With the notation of (2.5) one finds
\[
\Tr\big\{[\sL^{\dg{(o)}}_\psi] [\rK\dpsi]+ [\sL^{\dg{(o\star)}}_\psi] [\rK\dpsi]^\dagger +[\sL^{\dg{(\mu)}}_\psi][\rK\dpsi_{,\mu}]+
[\sL^{\dg{(\mu\star)}}_\psi] [\rK\dpsi_{,\mu}]^\dagger\big\}=0.
\]
If $\dpsi$ happens to be  a solution of the Lagrangian system, then with (2.16) this becomes
\[
\Tr\big\{[\frac{\pa}{\pa x^\mu}\sL^{\dg{(\mu)}}_\psi][\rK\dpsi]+ [\frac{\pa}{\pa x^\mu}\sL^{\dg{(\mu\star)}}_\psi][\rK\dpsi]^\dagger +[\sL^{\dg{(\mu)}}_\psi][\rK\dpsi]_{,\,\mu}+
[\sL^{\dg{(\mu\star)}}_\psi] [\rK\dpsi]^\dagger_{,\,\mu}\big\}=0,
\]
which leads to the wanted 'conserved current', since $\rK$ is supposedly constant.\hfill{$\blacksquare$}\\[3mm]
In gauge applications $\rK$ is often realized by a right multiplication by some $A\in\dC^{c\ti c}$. In such cases $\rK\dpsi$ in (4.13) should be replaced by $\dpsi A$.
\\[1mm]

All previous considerations can be applied to matrix gauge fields as well if we replace $\dpsi$ by $ \os{\aA} = \col[\ldots,\aA_\mu,\ldots]$. Some subtleties occur however because the range of
the functions $\aA_\mu$ is not the whole of $\dC^{c\ti c}$ but some real linear subspace $\gtg$ of it. See Appendix A for more details.\\[0.5mm]

This section is concluded with conservation laws for non-commutative free gauge fields which come from the special Lagrangian density (3.8).

\bt
 Consider the proto-Lagrangian $\sG$ of (3.6) with property (3.7) and Lagrange density as denoted in (3.8). For convenience restrict to $\gtg=\gtg^\da$ only.

 {\bf a.} Suppose $\sG^{(\dg{\nabla})}_A\!\cdot\!\osa=0$, for some $\osa\in\dR^N$ then we have the conservation law
 \be
 \sum_{\mu=1}^N\frac{\pa}{\pa x^\mu}\Big(\sum_{\ka=1}^N \Re\Tr\glh\proj\hat{\sG}^{(\dg{\mu\ka})}_A\!:\!(\osa\cdot\!\nabla)\aA_\ka\grh~-~a^\mu\sG_A\Big)\,=\,0\,.
 \ee

 {\bf b.} If for some $S=[S_\mu^\la]\in\dR^{N\ti N}$, with $\Tr S=0$, the assumptions
\be
\sG^{(\dg{\nabla})}_A\!\cdot \!S\osx=0~~~~\text{and}~~~~\Re\sum_{\mu,\,\nu =1}^N  \Tr\glh\,\hat{\sG}^{(\dg{\mu\nu})}_A\!:\!\sum_{\al=1}^N S_\mu^\al\pa_\al\aA_\nu\,\grh\,=\,0\,,
\ee
hold, then we have the conservation law
\be
\sum_{\mu=1}^N\frac{\pa}{\pa x^\mu}\Big(\sum_{\ka=1}^N 2\Re\Tr\glh\proj\hat{\sG}^{(\dg{\mu\ka})}_A (S\osx\cdot\nabla)\aA_\ka\grh\,-\,(S\osx\cdot\ose_\mu)\sG_A\Big)\,=\,0\,.
\ee
\et

{\bf Proof } \ \\
{\bf a.} Start from
\[
\gaf{}{s}\sG(~\ldots, \aF_{\mu\nu}(\osx-s\osa),\ldots;\ldots, \aF^\dagger_{\th\rho}(\osx-s\osa),\ldots;\osx-s\osa\,)\Big|_{s=0}=\gaf{}{s}\sG_A(\osx-s\osa)\Big|_{s=0}\,.
\]
Calculate the left hand side with the chain rule and use the assumptions
\[
-\sum_{\mu<\nu}\Tr\glh\sG^{(\dg{\mu\nu})}_A:(\osa\cdot\nabla)\aF_{\mu\nu}\grh~-~\sum_{\mu<\nu}\Tr\glh\sG^{(\dg{\mu\nu\star})}_A:(\osa\cdot\nabla)\aF_{\mu\nu}^\da\grh~-~\osa\cdot\sG^{\dg{\nabla}}_A~=
\]
\be
\hspace{5cm}=~-2\Re\sum_{\mu<\nu}\Tr\glh\sG^{(\dg{\mu\nu})}_A:(\osa\cdot\nabla)\aF_{\mu\nu}\grh\,.
\ee
With
\[
(\osa\cdot\nabla)\aF_{\mu\nu}=\pa_\mu(\osa\cdot\nabla\aA_\nu)-\pa_\nu(\osa\cdot\nabla\aA_\mu)-\llh\aA_\mu\lk\osa\cdot\nabla\aA_\nu\rlh+\llh\aA_\nu\lk\osa\cdot\nabla\aA_\mu\rlh\,,
\]
and the antisymmetries  $\mu\leftrightarrow\nu$, the expression (4.17) becomes, (mind the hat $ \hat{\mbox{}}~ $),
\[
-\Re\sum_{\mu,\nu=1}^N\Tr\glh\hat{\sG}^{(\dg{\mu\nu})}_A:\pa_\mu(\osa\cdot\nabla\aA_\nu)-\llh\aA_\mu\lk\osa\cdot\nabla\aA_\nu\rlh\grh\,=\hspace{5cm}
\]
\[
-\Re\sum_{\mu,\nu=1}^N\paf{}{x^\mu}\Tr\glh\hat{\sG}^{(\dg{\mu\nu})}_A\!:\!(\osa\cdot\nabla\aA_\nu)\grh~+~
\Re\sum_{\mu,\nu=1}^N\Tr\glh\pa_\mu\hat{\sG}^{(\dg{\mu\nu})}_A\!:\!(\osa\cdot\nabla\aA_\nu)~+~\hat{\sG}^{(\dg{\mu\nu})}_A\!:\!\llh\aA_\mu\lk\osa\cdot\nabla\aA_\nu\rlh\grh\,.
\]
The 2nd term is equal to
\[
\Re\sum_{\nu=1}^N\sum_{\mu=1}^N\Tr\glh\nabla_\mu^A\proj\hat{\sG}^{(\dg{\mu\nu})}_A\!:\!(\osa\cdot\nabla\aA_\nu)\grh\,=~0\,,
\]
because of the E-L-equations (3.13).

The right hand side of the 1st formula of this proof  equals $-\pa_\mu(a^\mu\sL_A)$. Hence (4.14).\\[1mm]
{\bf b.} Start from
\[
\gaf{}{s}\sG(~\ldots, \aF_{\mu\nu}(e^{sS}\osx),\ldots;\ldots, \aF^\dagger_{\th\rho}(e^{sS}\osx),\ldots;e^{sS}\osx\,)\Big|_{s=0}=\gaf{}{s}\sG_A(e^{sS}\osx)\Big|_{s=0}\,.
\]
Calculate the left hand side with the chain rule and use $\sG^{(\dg{\nabla})}_A\!\cdot \!S\osx=0$,
\[
2\Re\sum_{\mu<\nu}\Tr\glh\sG^{(\dg{\mu\nu})}_A:(S\osx\cdot\nabla)\aF_{\mu\nu}\grh\,=\,
\]
\[=\Re\sum_{\mu,\,\nu=1}^N\Tr\glh\hat{\sG}^{(\dg{\mu\nu})}_A:
\pa_\mu\big((S\osx\cdot\nabla)\aA_\nu\big)-\llh\aA_\mu\,\lk(S\osx\cdot\nabla)\aA_\nu\rlh-S^\al_\mu\pa_\al\aA_\nu\grh\,.
\]
Because of the assumption the very final contribution vanishes. Then we proceed as in part {\bf a}.\eindebewijs
\\[2mm]
{\bf Note} The orthogonality condition (4.15) is inspired by combining Thm 4.7 with Appendix A. Indeed, another way to obtain the preceding Theorem is to rewrite Thms 4.5, 4.7 in terms of $\os{\aA}$ with the aid of the table in Appendix A.

\bt
Consider the proto-Lagrangian $\sG$ of (3.6) with property (3.7) and Lagrange density as denoted in (3.8). For convenience consider $\gtg=\gtg^\da$ only. Suppose $\sG$ satisfies
\be
\sG(~\ldots, e^{s\aB}P_{\mu\nu}e^{-s\aB},\ldots;\ldots, e^{-s\aB^\dagger}P_{\th\rho}^\da e^{s\aB^\dagger},\ldots;\osx)=\sG(~\ldots, P_{\mu\nu},\ldots;\ldots, P_{\th\rho}^\da,\ldots;\osx)\,,
\ee
for all $P_{\mu\nu}\in\gtg\subset\dC^{c\ti c}\,,~1\leq\mu<\nu\leq N$,  some fixed $B\in\gtg$\, and (small) $s\in\dR$.

Then, for any solution $\osx\mapsto\ldots\aA_\mu(\osx)\ldots$ of the Lagrangian system of Theorem 3.3 one has the conservation law
\be
\sum_{\mu=1}^N\paf{}{x^\mu}\Re\Big(\sum_{\nu=1}^N\Tr\glh[\hat{\sG}^{(\dg{\mu\nu})}_A]:\llh\aB\lk\aA_\nu\rlh\grh\Big)~=~0\,.
\ee
\et
{\bf Proof~}
In (4.18) replace $P_{\mu\nu}\rightarrow\aF_{\mu\nu}$ and $Q_{\th\rho}\rightarrow\aF_{\th\rho}^\dagger$ and put the derivative to $s$ equal to $0$ at $s=0$,
\be
\sum_{1\leq\mu<\nu\leq N}\Tr\glh[\sG^{(\dg{\mu\nu})}_A]:(\aB\aF_{\mu\nu}-\aF_{\mu\nu}\aB)]\grh+\sum_{1\leq\th<\rho\leq N}\Tr\glh[\sG^{(\dg{\th\rho\star})}_A]:(-\aB^\dagger\aF_{\th\rho}^\dagger+\aF_{\th\rho}^\dagger\aB^\dagger)]\grh\,=0\,.
\ee
Due to the anti-symmetry in $\mu\leftrightarrow\nu$ of 
\[
\aB\aF_{\mu\nu}-\aF_{\mu\nu}\aB=\pa_\mu\llh\aB\lk\aA_\nu\rlh-\pa_\nu\llh\aB\lk\aA_\mu\rlh-\llh\aB\lk\llh\aA_\mu\lk\aA_\nu\rlh\rlh\,,
\]
applying convention (3.11), together with $\sG^{(\dg{\mu\nu}\star)}_A=[\sG^{(\dg{\mu\nu})}_A]^\dagger$, the 1st term of (4.20) equals the $\Re$-part of 
\[
\sum_{\mu=1}^N\sum_{\nu=1}^N\Tr\glh[\hat{\sG}^{(\dg{\mu\nu})}_A]:(\aB\aF_{\mu\nu}-\aF_{\mu\nu}\aB)\grh~=
\]

\[
=~\sum_{\mu=1}^N\sum_{\nu=1}^N\paf{}{x^\mu}\Tr\glh[\hat{\sG}^{(\dg{\mu\nu})}_A]\llh\aB\lk\aA_\nu\rlh\grh~-
~\sum_{\nu=1}^N\sum_{\mu=1}^N\paf{}{x^\nu}\Tr\glh[\hat{\sG}^{(\dg{\mu\nu})}_A]\llh\aB\lk\aA_\mu\rlh\grh~+~\,\hspace{1cm}
\]
\[
-\sum_{\nu=1}^N\sum_{\mu=1}^N\Tr\glh[\pa_\mu\hat{\sG}^{(\dg{\mu\nu})}_A]\llh\aB\lk\aA_\nu\rlh\grh~+
~\sum_{\mu=1}^N\sum_{\nu=1}^N\Tr\glh[\pa_\nu\hat{\sG}^{(\dg{\mu\nu})}_A]\llh\aB\lk\aA_\mu\rlh\grh~+
\]
\be
\hspace{2cm}~-~\sum_{\mu=1}^N\sum_{\nu=1}^N\Tr\glh[\hat{\sG}^{(\dg{\mu\nu})}_A]\llh\aB\lk\llh\aA_\mu\lk\aA_\nu\rlh\rlh\grh~.
\ee
On the 2nd line we apply the E-L-equations (3.13) together with $\pa_\nu\hat{\sG}^{(\dg{\mu\nu})}_A=-\pa_\nu\hat{\sG}^{(\dg{\nu\mu})}_A\,$. This together with the 3rd line leads to
\[
-\sum_{\nu=1}^N\sum_{\mu=1}^N\Tr\glh\llh\aA_\mu\lk\hat{\sG}^{(\dg{\mu\nu})}_A\rlh\llh\aB\lk\aA_\nu\rlh\grh~+
~\sum_{\mu=1}^N\sum_{\nu=1}^N\Tr\glh\llh\aA_\nu\lk\hat{\sG}^{(\dg{\mu\nu})}_A\rlh\llh\aB\lk\aA_\mu\rlh\grh~+
\]
\[
\hspace{2cm}~-~\sum_{\mu=1}^N\sum_{\nu=1}^N\Tr\glh[\hat{\sG}^{(\dg{\mu\nu})}_A]\llh\aB\lk\llh\aA_\mu\lk\aA_\nu\rlh\rlh\grh~.
\]
These 3 terms add up to $0$ because for each pair $\mu,\nu$ separately we can apply the  identity
\be
-\Tr\glh\llh\mM\lk\mG\rlh:\llh\mB\lk\mN\rlh\grh\,+\,\Tr\glh\llh\mN\lk\mG\rlh:\llh\mB\lk\mM\rlh\grh\,=\,\Tr\glh\mG:\llh\mB\lk\llh\mM\lk\mN\rlh\rlh\grh\,,
\ee
for matrices $\mG,\mB,\mM,\mN\in\dC^{r\ti r}$\,.

(Of course the two terms on the 3rd line of (4.21) are equal. But then, using that equality,  the latter trick no longer works for each index pair $\mu,\nu$ separately!)

Thus we found out that  (4.20) corresponds to (4.19). \eindebewijs

\section{Static/Dynamic Gauge Extensions of Lagrangians}

A basic ingredient for this section is a (fixed) Lie-group $\gtG\subset\dC^{c\ti c}$ of invertible $c\ti c$-matrices. Its Lie-algebra $\gtg$ is a $\dR$-linear subspace of $\dC^{c\ti c}$.
Important examples are (subgroups of) $\gtG_J$, for some fixed invertible matrix $J\in\dC^{c\ti c}$. The relevant definitions are as in section 3,
\be
\gtG_J\,=\,\big\{~\gU\in\dC^{c\ti c}~\big|~\gU^\da J\gU=J~\big\}\,,~~~~~\gtg_J\,=\,\big\{~\aA\in\dC^{c\ti c}~\big|~\aA^\da J +J\aA =0~\big\}\,.
\ee

In the discussion to follow suitable subspaces of

\[ \text{the group} ~~~\gtG_{\text{loc}}=\fG~~~ ~~~~~~\text{and the $\dR$-linear space}~~~~\fg
\]

will be used.
It will be tacitly assumed that the behaviour at $\infty$ of the considered subspaces
is such that our formulae make sense. The $\ms{C}^\infty$-smoothness condition can often be relaxed. Neither of those assumptions will bother us.

The group action from the right of $\fG$  on $\fpsi$ is naturally defined by
\[
\fpsi\ti\fG~\to~\fpsi~:~~~(\dpsi \cU)(\osx)=\dpsi(\osx)\cU(\osx).
\]
For each $1\leq\mu\leq N$, a group action from the right of $\fG$  on $\fg$ is  defined by
\[
\fg\ti\fG~\to~\fg~:~~(\cA_\mu\dhr \cU)(\osx)=\cU^{-1}(\osx)\cA_\mu(\osx) \cU(\osx)-\cU^{-1}(\osx)(\pa_\mu \cU)(\osx)\,.
\]
In the proof of Thm 1.2 it has been shown that this action ('gauge transform')is indeed a (inhomogeneous) group action. This means
\be
[\cA_\mu\dhr \cU]\dhr \cV\,=\,\cA_\mu\dhr(\cU\cV)~~~~\,.
\ee
As before, for given $\cA_\mu,\cA_\nu\in\fg\,,1\leq\mu,\nu\leq N$, define
\be
\cF_{\mu\nu}=\pa_\mu\cA_\nu-\pa_\nu\cA_\mu-\llh\cA_\mu\lk\cA_\nu\rlh~\in~\fg\,.
\ee
Then
\be
\cU^{-1}\cF_{\mu\nu}\cU=\pa_\mu(\cA_\nu\dhr \cU)-\pa_\nu(\cA_\mu\dhr \cU)-\llh(\cA_\mu\dhr \cU)\lk(\cA_\nu\dhr \cU)\rlh\,.
\ee

\bt
 Fix a matrix Lie-Group $\gtG\subset\dC^{c\ti c}$. Suppose a proto-Lagrangian $\sL$, cf. (2.13), to be $\gtG$-invariant, i.e.
\footnote{Property (5.5) is named {\em Global Gauge Invariance} by physicists. The conclusion of Theorem 5.1 is named, in physicists' vernacular, the property of {\em Local Gauge Invariance}.
In mathematicians' jargon however, the usage of 'global', as opposed to 'local', usually refers to a more involved (more difficult) notion.}
\[
\forall\,\gU\in\gtG~~\forall\,\mP\in\dC^{r\ti c}~~\forall\,\os{\mR}\in\dC^{Nr\ti c}~~\forall\,\osx\in\dR^N\,: \hspace{7cm}
\]
\be
\sL(\mP \gU\,;\,\gU^\da\mP^\da\,;\,\os{\mR}\gU\,;\,\gU^\da\os{\mR}^\da\,;\,\osx)=\sL(\mP\,;\,\mP^\da\,;\,\os{\mR}\,;\,\os{\mR}^\da\,;\,\osx)\,
\ee
Then, for all $\osx\in\dR^N$, the {\bf statically gauge extended} Lagrangian density
\be
\sL_{\psi,\,A}(\osx)=\sL(\dpsi\,;\dpsi^\da\,;\ldots,\pa_\mu\dpsi+\dpsi\cA_\mu\,,\ldots\,;\ldots,\pa_\mu\dpsi^\da+\cA_\mu^\da\dpsi^\da\,,\ldots\,;\osx)\,,\hspace{2cm}
\ee
\[
 \hspace{3cm}\text{with any}~~\dpsi\in\fpsi\,,~\cA_\mu\in\fg\,,1\leq\mu\leq N\,,
\]
{\bf equals} the {\bf statically gauge extended} Lagrangian density
\be
\sL_{\psi\cU,\,A\dhr\,\,\cU}(\osx)= \hspace{10cm}
\ee
\[
=\sL(\dpsi\cU\,;\cU^\da\dpsi^\da\,;\ldots,\pa_\mu(\dpsi\cU)+(\dpsi\cU)(\cA_\mu\dhr\cU)\,,\ldots\,;\ldots,\pa_\mu(\dpsi\cU)^\da+
(\cA_\mu\dhr\cU)^\da(\dpsi\cU)^\da,\ldots\,;\osx)\,,
\]
\[
\hspace{-2cm}\text{with any}~~~~~\cU\in\fG\,.
\]
In (5.6),(5.7) we wrote $\dpsi$ instead of $\dpsi(\osx)$, etc.
\et
{\bf Proof } Straightforward calculation. \eindebewijs


\begin{examp}
Consider the proto-Lagrangian, cf. (2.13), 
\[
\sL(\mP;\mQ^\top;\os{\mR}\,;\os{\mS}^\top;\,\osx)\,=\, \ri\,\Tr[\mQ^\top(\sum_\mu\Ga^\mu\mR_\mu+M\mP)]\,
\]
with fixed $\Ga^\mu,M\in\dC^{r\ti r}$ and $[\Ga^\mu]^\da=\Ga^\mu\,,M^\da=-M$. Put $\gtG=\gt{U}(c)\subset\dC^{c\ti c}$, that is the unitary group $\gtG_I$, with $I$ the identity matrix.
Our proto-Lagrangian is $\gt{U}(c)$-invariant
\[
 \ri\,\Tr[\gU^\da\mP^\da(\Ga^\mu\mR_\mu\gU+M\mP\gU)]\,=\, \ri\,\Tr[\mP^\da(\Ga^\mu\mR_\mu+M\mP)]\,,~~~\gU\in\gt{U}(c)\,,
\]
because $\gU^\da=\gU^{-1}$ and the properties of $\Tr$.

Then the statically extended Lagrangian density
\be
\sL_{\psi,\,A}(\osx)= \ri\,\Tr[\dpsi^\da\big(\Ga^\mu(\pa_\mu\dpsi+\dpsi\aA_\mu)+M\dpsi\big)]\,,\hspace{4cm}
\ee
\[
 \hspace{3cm}\text{with any}~~\dpsi\in\fpsi\,,~\cA_\mu\in\ms{C}^\infty(\dR^N\!:\!\gt{u}(c))\,,1\leq\mu\leq N\,,
\]
 equals the statically extended Lagrangian density
\be
\sL_{\psi\cU,\,A\dhr\,\cU}(\osx)= \ri\, \Tr[\cU^\da\dpsi^\da\big(\Ga^\mu(\pa_\mu(\dpsi\cU)+\dpsi\cU(\cU^{-1}\aA_\mu\cU-\cU^{-1}\pa_\mu\cU))+M\dpsi\cU\big)]\,,
\ee
\[
\hspace{6cm}\text{with any}~~~~~\cU\in\ms{C}^\infty(\dR^N\!:\!\gt{U}(c))\,.
\]
Note that, if $M$ is replaced by the 'nonlinearity' $\ri\dpsi\dpsi^\da$, the argument still holds. \eindebewijs
\end{examp}

\bt
$\bullet~$ Suppose that the statically gauge extended Lagrange density $\sL_{\psi,\,A}$, cf. (5.6) leads to an $\dR$-valued Langrangian functional $\cL_{\psi,\,A}$. The  E-L-equations are
\be
\begin{array}{c} \displaystyle
\sL^{\dg{(o)}}_{\psi,A}-\sum_{\mu=1}^N \Big(\frac{\pa}{\pa x^\mu}[\sL^{\dg{(\mu)}}_{\psi,A}]-[\aA_\mu\sL^{\dg{(\mu)}}_{\psi,A}]\Big)\,=0\,,\\[4mm]
\proj\Big(\dpsi^\da[\sL^{\dg{(\ka)}\da}_{\psi,A}+\sL^{\dg{(\ka\star)}}_{\psi,A}]\Big)=0~,~~~~~~
\proj\Big(\dfrac{\dpsi^\da[\sL^{\dg{(\ka)}\da}_{\psi,A}-\sL^{\dg{(\ka\star)}}_{\psi,A}]}{\ri}\Big)=0~,~~~1\leq\ka\leq N\,.
\end{array}
\ee
Here $\proj: \dC^{c\ti c}\to\dC^{c\ti c}$ denotes the $\dR$-orthogonal projection on $\gtg$.

$\bullet$ If it happens that $\proj(\ri\mZ)=\ri\proj^\perp\mZ\,,\,\mZ\in\dC^{c\ti c}$, the 2nd line in (5.10) reduces to
\be
\dpsi^\da\sL^{\dg{(\ka)}\da}_{\psi,A}+(\proj-\proj^\perp)\dpsi^\da\sL^{\dg{(\ka\star)}}_{\psi,A}=0~,~~~1\leq\ka\leq N\,.
\ee

$\bullet$ In the important special case $\gtg=\gtg_J$, with $J=J^\da=J^{-1}$, (5.11) can be written
\be
\sL^{\dg{(\ka)}}_{\psi,A}\dpsi-J\dpsi^\da\sL^{\dg{(\ka\star)}}_{\psi,A}J\,=\,0~,~~~1\leq\ka\leq N\,.
\ee
\et
{\bf Proof } $\bullet~$ The perturbed statically extended Lagrangian  $\sL_{\psi,\,A}$ reads
\[
\sL(\dpsi +\vep\dk{H}\,;\dpsi^\da+\vep^\star\dk{K}\,;\ldots,\pa_\mu(\dpsi +\vep\dk{H})+(\dpsi+\vep\dk{H})(\cA_\mu+\vep_\ka\del_{\mu\ka}\cH) \,,\ldots\,; \hspace{2cm}
\]
\[\hspace{5cm};\ldots,\pa_\mu(\dpsi^\da+\vep^\star\dk{K})+(\cA_\mu^\da+\vep_\ka\del_{\mu\ka}\cH^\da)(\dpsi^\da+\vep^\star\dk{K})\,,\ldots\,;\osx)\,
\]
The results of $\gaf{}{\vep}\big|_{\vep=0}\,,~~\gaf{}{\vep^\star}\big|_{\vep^\star=0}\,,~~\gaf{}{\vep_\ka}\big|_{\vep_\ka=0}\,,~1\leq\ka\leq N$, being put to $0$ are,

for all functions $\dk{H}\,,\dk{K}\,,\aH$\,,
 \[
\Tr\big[\sL^{\dg{(o)}}\!:\!\dk{H}\big]\,+\,\sum_\mu\Tr\big[\sL^{\dg{(\mu)}}\!:\!\pa_\mu\dk{H}\big]\,+\,\sum_\mu\Tr\big[\sL^{\dg{(\mu)}}\!:\!\dk{H}\aA_\mu\big]\,=\,0\,,
\]
\[
\Tr\big[\sL^{\dg{(o\star)}}\!:\!\dk{K}\big]\,+\,\sum_\mu\Tr\big[\sL^{\dg{(\mu\star)}}\!:\!\pa_\mu\dk{K}\big]\,+\,\sum_\mu\Tr\big[\sL^{\dg{(\mu\star)}}\!:\!\aA_\mu^\da\dk{K}\big]\,=\,0\,,
\]
\[
\sum_\mu\Tr\big[\sL^{\dg{(\mu)}}\!:\!\dpsi\del_{\mu\ka}\aH\big]\,+\,\sum_\mu\Tr\big[\sL^{\dg{(\mu\star)}}\!:\!\del_{\mu\ka}\aH^\da\dpsi^\da\big]\,=\,0\,, ~1\leq\ka\leq N\,.
\]
The usual partial integration techniques applied to the first two lines lead to the E-L-equations for $\dpsi$. Also use Theorem 2.4.

From the final line we arrive at (5.10) because of the trace identity
\be
\Tr\glh\mX\mZ+\mY\mZ^\da\grh\,=\,\Re\Tr\glh\Big(\mX^\da+\mY\Big)^\da\mZ\grh\,-\,\ri\,\Re\Tr\glh\Big(\frac{\mX^\da-\mY}{\ri}\Big)^\da\mZ\grh\,.
\ee
$\bullet~$ If for $\mX,\mY\in\dC^{c\ti c}$ one has $\proj(\mX+\mY)=0$ and $\proj^\perp(\mX-\mY)=0$, it follows that $\mX+(\proj-\proj^\perp)Y=0$ and also $\mY+(\proj-\proj^\perp)X=0$.

$\bullet~$ In this special case $(\proj-\proj^\perp)Y=-J\mY^\da J$ and $\proj[\mY^\da]=[\proj\mY]^\da$.
\eindebewijs

\  \\[4mm]
\begin{examps} \ \\
Note that in the E-L-equations (5.10) the $\cA_\mu$ occur only 'algebraically'.\\ The $\pa_\mu\cA$ are not involved!\\[0.5mm]
{\bf a.} For the Lagrangian densities from examples 2.5a and 5.2 the 2nd set of E-L-equations (5.12) does not depend on $\cA$. If we choose $\gtg=\gtg_J$, the 2nd line reads
\[
\dpsi^\da\Ga^\ka\dpsi\,=\,0\,,~~~~1\leq\ka\leq N\,.
\]
It means that $\dpsi$ can only take values in a cone in $\dC^{r\ti c}$. If one of the $\Ga^\ka=\Ga^{\ka\da}$ is strictly positive, the only solutions are $\dpsi=0$, the trivial ones.
If a nontrivial choice for $\dpsi$ is possible it can be substituted in the 1st E-L-equation and we are left with  an algebraic equation for the $\cA_\ka$. \\[0.5mm]
{\bf b.}  For the Lagrangian densities from example 2.5c, again with $\gtg=\gtg_J~$, the 2nd set of E-L-equations becomes
\[
\sum_{\mu=1}^N[\pa_\mu\dpsi+\dpsi\cA_\mu]^\da\Th^{\mu\ka}\dpsi\,-\,J\Big( \sum_{\mu=1}^N[\dpsi^\da\Th^{\ka\mu}[\pa_\mu\dpsi+\dpsi\cA_\mu]\Big)J\,=\,0\,,~~~~1\leq\ka\leq N\,,
\]
which is algebraic in the $\cA_\ka$\,.  \eindebewijs
\end{examps}

\ \\[3mm]

Finally we want to consider the {\bf dynamically gauge extended} Lagrangian density or

{\bf Gauge field extended} Lagrangian density of type $\sL_{\psi,\,A}(\osx)+\sG_A(\osx)\,$.

\bt
Fix a matrix Liegroup $\gtG\subset\dC^{c\ti c}$ with Lie algebra $\gtg\subset\dC^{c\ti c}$ and property $\gtg^\da=\gtg$.

Fix  a proto Lagrangian of type (2.13)
\[(\mP;\mQ^\top;\os{\mR}\,;\os{\mS}^\top;\,\osx)~~\mapsto~~\sL(\mP;\mQ^\top;\os{\mR}\,;\os{\mS}^\top;\,\osx)\,,
\]
 leading to a $\dR$-valued Lagrangian functional $\cL$. Require the special property
\be
\forall\,\mP~\forall\,\os{\mR}~\forall\,\osx:~ \proj\Big(\frac{\mP^\da\big[\sL^{\dg{(\ka)}\da}(\mP;\mP^\da;\os{\mR}\,;\os{\mR}^\da;\,\osx)\,-\,\sL^{\dg{(\ka\star)}}(\mP;\mP^\da;\os{\mR}\,;
\os{\mR}^\da;\,\osx)\big]}{\ri}\Big)~=~0\,.
\ee

Fix a second proto Lagrangian of type (3.6) and such that
\[
\forall\,R_{\mu\nu}\in\gtg~:~~\sG(\ldots,R_{\mu\nu},\ldots;\ldots,R_{\th\rho}^\dagger,\ldots;\osx)\,\in\,\dR\,.
\]
Consider the {\bf dynamically extended} Lagrangian density
\[
\sL_{\psi,\,A}(\osx)+\sG_A(\osx)\,=\,\sL(\dpsi\,;\dpsi^\da\,;\ldots,\pa_\mu\dpsi+\dpsi\cA_\mu\,,\ldots\,;\ldots,\pa_\mu\dpsi^\da+\cA_\mu^\da\dpsi^\da\,,\ldots\,;\osx)\,+\hspace{1cm}
\]
\be
\hspace{7cm}+\,\sG(~\ldots, \aF_{\mu\nu}(\osx),\ldots;\ldots, \aF^\dagger_{\th\rho}(\osx),\ldots;\osx\,)
\ee
\[
 \hspace{0cm}\text{with any}~~~\dpsi\in\fpsi\,,~~\cA_\mu\in\fg\,,~1\leq\mu\leq N\,.\hspace{7cm}
\]

$\bullet~$ The Euler-Lagrange equations are, with $\sL^{\dg{(o)}}_{\psi,A}$ instead of $\sL^{\dg{(o)}}_{\psi,A}(\osx)$, etc.,
\be
\begin{array}{c} \displaystyle
[\sL^{\dg{(o)}}_{\psi,A}]-\sum_{\mu=1}^N \Big(\frac{\pa}{\pa x^\mu}[\sL^{\dg{(\mu)}}_{\psi,A}]-[\aA_\mu\sL^{\dg{(\mu)}}_{\psi,A}]\Big)\,=0\,,\\[5mm]
\proj\Big(\dpsi^\da[\sL^{\dg{(\ka)}\da}_{\psi,A}+\sL^{\dg{(\ka\star)}}_{\psi,A}]\Big)\,-\,2\sum_{\mu=1}^N\Big(\pa_\mu\proj[\hat{\sG}^{(\dg{\mu\ka})}_A]-
\llh\cA_\mu\lk\proj[\hat{\sG}^{(\dg{\mu\ka})}_A]\rlh\Big)^\da~=~0\,,~~~1\leq\ka\leq N\,.\\
\end{array}
\ee
Here $\proj: \dC^{c\ti c}\to\dC^{c\ti c}$ denotes the $\dR$-orthogonal projection on $\gtg$.\\[1mm]
$\bullet~$ In the  special case $\gtg=\gtg_J$, with $J=J^\da=J^{-1}$, the 2nd line in (5.16) can be rewritten
\be
\sL^{\dg{(\ka)}}_{\psi,A}\dpsi-J\dpsi^\da\sL^{\dg{(\ka\star)}}_{\psi,A}J\,\,-\,2\sum_{\mu=1}^N\Big(\pa_\mu\proj[\hat{\sG}^{(\dg{\mu\ka})}_A]-
\llh\cA_\mu\lk\proj[\hat{\sG}^{(\dg{\mu\ka})}_A]\rlh\Big)~=~0\,,~~~1\leq\ka\leq N\,.
\ee
\et
{\bf Proof }
$\bullet~$ The perturbed gauge supplemented Lagrangian reads
\[
\sL(\dpsi +\vep\dk{H}\,;\dpsi^\da+\vep^\star\dk{K}\,;\ldots,\pa_\mu(\dpsi +\vep\dk{H})+(\dpsi+\vep\dk{H})(\cA_\mu+\vep_\ka\del_{\mu\ka}\cH) \,,\ldots\,; \hspace{6.5cm}
\]
\[ \hspace{5cm} ;\ldots,\pa_\mu(\dpsi^\da+\vep^\star\dk{K})+(\cA_\mu^\da+\vep_\ka\del_{\mu\ka}\cH^\da)(\dpsi^\da+\vep^\star\dk{K})\,,\ldots\,;\osx)\,+ \hspace{2.5cm}
\]
\[
\hspace{7cm}+~\sG(\ldots,\aF_{\mu\nu,\vep\ka}\,,\ldots;\ldots,\aF_{\th\rho,\vep\ka}^\dagger,\ldots;\osx)\,,~~~~1\leq\ka\leq N\,,
\]
where
\[\aF_{\mu\nu;\vep,\ka}=\aF_{\mu\nu}+\vep_\ka\,\del_{\nu\ka}\glh\pa_\mu\cH-\llh\cA_\mu \lk \cH\rlh\,\grh\,
-\vep_\ka\,\del_{\mu\ka}\glh\pa_\nu\cH-\llh\cA_\nu \lk \cH\rlh\,\grh\,,
\]
The results of $\gaf{}{\vep}\big|_{\vep=0}\,,~~\gaf{}{\vep^\star}\big|_{\vep^\star=0}\,~~\gaf{}{\vep_\ka}\big|_{\vep_\ka=0}$, being put to $0$ are, respectively,
\[
\Tr\big[\sL^{\dg{(o)}}\!:\!\dk{H}\big]\,+\,\sum_\mu\Tr\big[\sL^{\dg{(\mu)}}\!:\!\pa_\mu\dk{H}\big]\,+\,\sum_\mu\Tr\big[\sL^{\dg{(\mu)}}\!:\!\dk{H}\aA_\mu\big]\,=\,0\,,\hspace{3.5cm}
\]
\[
\Tr\big[\sL^{\dg{(o\star)}}\!:\!\dk{K}\big]\,+\,\sum_\mu\Tr\big[\sL^{\dg{(\mu\star)}}\!:\!\pa_\mu\dk{K}\big]\,+\,\sum_\mu\Tr\big[\sL^{\dg{(\mu\star)}}\!:\!\aA_\mu^\da\dk{K}\big]\,=\,0\,,\hspace{3cm}
\]
\[
\sum_\mu\Tr\big[\sL^{\dg{(\mu)}}\!:\!\dpsi\del_{\mu\ka}\aH\big]\,+\,\sum_\mu\Tr\big[\sL^{\dg{(\mu\star)}}\!:\!\del_{\mu\ka}\aH^\da\dpsi^\da\big]\,+\,\hspace{5.5cm}
\]
\[
\hspace{1cm}-~2\sum_{\mu}\Re\Tr\glh\Big(\proj\pa_\mu\hat{\sG}^{(\dg{\mu\ka\star})}_A+\proj\llh\cA^\dagger_\mu\,\lk\proj\hat{\sG}^{(\dg{\mu\ka\star})}_A\rlh\Big)^\dagger[\aH]\grh\,=\,0\,,~~1\leq\ka\leq N\,.
\]
With (5.13) the 3rd set of  equations can be rewritten
\[
\Re\Tr\glh\big(\dpsi^\da([\sL^{\dg{(\ka)}}]^\da+[\sL^{\dg{(\ka\star)}}]\big)^\da\aH\grh\,+\,\ri\Re\Tr\glh\big(\ri\dpsi^\da([\sL^{\dg{(\ka)}}]^\da-[\sL^{\dg{(\ka\star)}}]\big)^\da\aH\grh\,+\,
\]
\[
\hspace{1cm}-~2\sum_{\mu}\Re\Tr\glh\Big(\proj\pa_\mu\hat{\sG}^{(\dg{\mu\ka})}_A\,-\,\llh\cA_\mu\,\lk\proj\hat{\sG}^{(\dg{\mu\ka})}_A\rlh\Big)^{\da\,\da}[\aH]\grh\,=\,0\,,~~1\leq\ka\leq N\,.
\]
Because of assumption (5.14) the $\ri\Re\Tr$-term cancels. The assumption $\gtg^\da=\gtg$ enables us to interchange $\da$ and $\proj$.

$\bullet~$ Finally (5.17) follows as in the proof of Thm (5.3). \eindebewijs

\ \\[2mm]

Finally we want to find the conservation law of 'conserved currents'.

\bt
Consider proto-Lagrangians $\sL$ and $\sG$ as in Theorem 5.5. Suppose for some $\aB\in\gtg$ they {\em both} have the invariance properties
\[ \sL(\mP e^{s\aB};(\mP e^{s\aB})^\dagger;\ldots \mQ_\la e^{s\aB}\ldots;\ldots(\mQ_\la e^{s\aB})^\dagger\ldots;\,\osx) =\hspace{3cm}
\]
\be
~\hspace{5cm}=\sL(\mP;\mP^\dagger;\ldots \mQ_\la\ldots;\ldots\mQ_\la^\dagger\ldots;\,\osx)+\ms{O}(s^2)\,,
\ee

\[
\sG(~\ldots, e^{-s\aB}\mR_{\mu\nu}e^{s\aB},\ldots;\ldots, e^{s\aB^\da}\mR^\dagger_{\th\rho}e^{-s\aB^\da},\ldots;\osx\,)\,=\,\hspace{3cm}
\]
\be
\hspace{6cm }=\,\sG(~\ldots,\mR_{\mu\nu},\ldots;\ldots, \mR^\dagger_{\th\rho},\ldots;\osx\,)\,+\ms{O}(s^2)\,.
\ee
Then, the solutions to the E-L-system (5.16) satisfy the conservation law
\be
\sum_{\mu=1}^N\paf{}{x_\mu}\Big\{\,\Tr\glh\sL^{\dg{(\mu)}}_{\psi,A}\!:\!\dpsi\aB\grh\,+\,\Tr\glh\sL^{\dg{(\mu\star)}}_{\psi,A}\!:\!\aB^\da\dpsi^\da\grh+
\sum_{\ka=1}^N 2\Re\Tr\glh\proj\hat{\sG}^{(\dg{\mu\ka})}_A:\llh\aA_\ka\lk\aB\rlh\grh\,\Big\}~=~0.
\ee
\et

{\bf Proof~}
Add the Lagrange densities $\sL_{\psi,A}$ and $\sG_A$  and  put to $0$ the $\gaf{}{s}$ of the expression
\[
\sL(\dpsi e^{s\aB}\,;e^{s\aB^\da}\dpsi^\da\,;\ldots,\pa_\mu\dpsi e^{s\aB}+\dpsi\cA_\mu e^{s\aB}\,,\ldots\,;\ldots,e^{s\aB^\da}\pa_\mu\dpsi^\da+e^{s\aB^\da}\cA_\mu^\da\dpsi^\da\,,\ldots\,;\osx)\,+\hspace{1cm}
\]
\[
\hspace{6.5cm} +\,\sG(~\ldots, e^{-s\aB}\aF_{\mu\nu}e^{s\aB},\ldots;\ldots, e^{s\aB^\da}\aF^\dagger_{\th\rho}e^{-s\aB^\da},\ldots;\osx\,)\,
\]
One finds, 

\[
\Tr\Big[\sL^{\dg{(o)}}_{\psi,A}\!:\!\dpsi\aB\grh+\sum_\mu\Tr\glh\sL^{\dg{(\mu)}}_{\psi,A}\!:\!\pa_\mu\!\dpsi\aB\grh+\sum_\mu\Tr\glh\sL^{\dg{(\mu)}}_{\psi,A}\!:\!\dpsi\aA_\mu\aB~\grh~+ \hspace{5cm}
\]
\[+~ \Tr\glh\sL^{\dg{(o\star)}}_{\psi,A}\!:\!\aB^\da\dpsi^\da\grh
+\sum_\mu\Tr\glh\sL^{\dg{(\mu\star)}}_{\psi,A}\!:\!\aB^\da\pa_\mu\!\dpsi^\da\grh+\sum_\mu\Tr\glh\sL^{\dg{(\mu\star)}}_{\psi,A}\!:\!\aB^\da\aA_\mu^\da\dpsi^\da\grh~+~\hspace{1.5cm}
\]
\be\hspace{3cm}~+~\sum_{\mu<\nu}\Tr\glh\sG^{(\dg{\mu\nu})}_A\!:\!\llh\aF_{\mu\nu}\lk\aB\rlh\grh~+~
\sum_{\th<\rho}\Tr\glh\sG^{(\dg{\th\rho\star})}_A\!:\!\llh\aB^\da\lk\aF_{\th\rho}^\da\rlh\grh~=~0\,.
\ee
Rewrite the 3rd term and the 6th term:
\[
\sum_\mu\Tr\glh\sL^{\dg{(\mu)}}_{\psi,A}\!:\!\dpsi\aA_\mu\aB\grh\,=\,\sum_\ka\Tr\glh\sL^{\dg{(\ka)}}_{\psi,A}\!:\!\dpsi\llh\aA_\ka\lk\aB\rlh\grh\,+
\sum_\mu\Tr\glh\aA_\mu\sL^{\dg{(\mu)}}_{\psi,A}\!:\!\dpsi\aB\grh\,,
\]
\[
\sum_\mu\Tr\glh\sL^{\dg{(\mu\star)}}_{\psi,A}\!:\!(\dpsi\aA_\mu\aB)^\da\grh\,=\,\sum_\ka\Tr\glh\sL^{\dg{(\ka\star)}}_{\psi,A}\!:\!(\dpsi\llh\aA_\ka\lk\aB\rlh)^\da\grh\,+
\sum_\mu\Tr\glh\aA_\mu^\da\sL^{\dg{(\mu\star)}}_{\psi,A}\!:\!(\dpsi\aB)^\da\grh\,.
\]
These identities, together with the 1st E-L-equation of (5.16) turn the first 6 terms of (5.21) into
\[
\sum_\mu\pa_\mu\Tr\glh\sL^{\dg{(\mu)}}_{\psi,A}\!:\!\dpsi\aB\grh + \sum_\mu\pa_\mu\Tr\glh\sL^{\dg{(\mu\star)}}_{\psi,A}\!:\!\aB^\da\dpsi^\da\grh + \hspace{5cm}
\]
\[
\hspace{3cm} +\sum_\ka\Tr\glh\sL^{\dg{(\ka)}}_{\psi,A}\!:\!\dpsi\llh\aA_\ka\lk\aB\rlh\grh\,+\sum_\ka\Tr\glh\sL^{\dg{(\ka\star)}}_{\psi,A}\!:\!(\dpsi\llh\aA_\ka\lk\aB\rlh)^\da\grh\,
\]
With Trace identity (5.13) and condition (5.14) the latter becomes
\[
\sum_\mu\pa_\mu\Tr\glh\sL^{\dg{(\mu)}}_{\psi,A}\!:\!\dpsi\aB\grh + \sum_\mu\pa_\mu\Tr\glh\sL^{\dg{(\mu\star)}}_{\psi,A}\!:\!\aB^\da\dpsi^\da\grh + \hspace{5cm}
\]
\be
\hspace{3cm}+~2\sum_{\ka,\,\mu=1}^N\Re\Tr\glh\Big(\proj\pa_\mu\hat{\sG}^{(\dg{\mu\ka})}_A\,-\,\llh\cA_\mu\,\lk\proj\hat{\sG}^{(\dg{\mu\ka})}_A\rlh\Big):\llh\aA_\ka\lk\aB\rlh\grh\,.
\ee
Next, because of (anti)symmetry,   $\aB\in\gtg$ being constant and the definition of $\cF_{\mu\nu}$, the final 2 terms of (5.21)  equal to 

\[ \Re\sum_{\mu,\nu=1}^N\Tr\glh\hat{\sG}^{(\dg{\mu\nu})}_A\!:\!\llh\aF_{\mu\nu}\lk\aB\rlh\grh~=~
\Re\sum_{\mu,\nu=1}^N\Tr\glh\hat{\sG}^{(\dg{\mu\nu})}_A:\pa_\mu\llh\aA_\nu\lk\aB\rlh\grh~+ \hspace{3cm}
\]
\[
\hspace{1.5cm}-~\Re\sum_{\mu,\,\nu=1}^N\Tr\glh\hat{\sG}^{(\dg{\mu\nu})}_A\!\!:\!\pa_\nu\llh\aA_\mu\lk\aB\rlh\grh
~-~\Re\sum_{\mu,\,\nu=1}^N\Tr\glh\sG^{(\dg{\mu\nu})}_A\!:\!\llh\llh\aA_\mu\lk\aA_\nu\rlh\lk\aB\rlh\grh\,=
\]
\be
=~2\Re\sum_{\mu,\nu=1}^N\Tr\glh\hat{\sG}^{(\dg{\mu\nu})}_A:\pa_\mu\llh\aA_\nu\lk\aB\rlh\grh\,
-~\Re\sum_{\mu,\,\nu=1}^N\Tr\glh\sG^{(\dg{\mu\nu})}_A\!:\!\llh\llh\aA_\mu\lk\aA_\nu\rlh\lk\aB\rlh\grh\,.\hspace{1.5cm}
\ee
If we add (5.22), (5.23), we arrive at (5.20), up to a term
\[
-~\Re\,\sum_{\ka,\,\mu=1}^N\Big(\,2\,\Tr\glh\llh\cA_\mu\,\lk\proj\hat{\sG}^{(\dg{\mu\ka})}_A\rlh:\llh\aA_\ka\lk\aB\rlh\grh\,+
\Tr\glh\proj\sG^{(\dg{\mu\ka})}_A\!:\!\llh\llh\aA_\mu\lk\aA_\ka\rlh\lk\aB\rlh\grh\,\Big)\,.
\]
Split the first term in this summation. It  becomes,
\[
-~\Re\,\sum_{\ka,\,\mu=1}^N\Big(\,\Tr\glh\llh\cA_\mu\,\lk\proj\hat{\sG}^{(\dg{\mu\ka})}_A\rlh:\llh\aA_\ka\lk\aB\rlh\grh\,
-\Tr\glh\llh\cA_\ka\,\lk\proj\hat{\sG}^{(\dg{\mu\ka})}_A\rlh:\llh\aA_\mu\lk\aB\rlh\grh\,+ \hspace{1.5cm}
\]
\[\hspace{9cm} +~
\Tr\glh\proj\sG^{(\dg{\mu\ka})}_A\!:\!\llh\llh\aA_\mu\lk\aA_\ka\rlh\lk\aB\rlh\grh\,\Big)\,.
\]
Each term in this sum equals $0$ because of the trace identity
\[
\Tr\glh\llh\mM\lk\mG\rlh:\llh\mK\lk\mB\rlh\grh\,-\,\Tr\glh\llh\mK\lk\mG\rlh:\llh\mM\lk\mB\rlh\grh\,+\,\Tr\glh\mG:\llh\llh\mM\lk\mK\rlh\lk\mB\rlh\grh\,=\,0.
\]
Indeed, note that for any $\mM,\mG,\mK,\mB\,\in\dC^{c\ti c}$,
\[
\Tr\glh~\mM\mG\mK\mB-\mG\mM\mK\mB-\mM\mG\mB\mK+\mG\mM\mB\mK-\mK\mG\mM\mB+\mG\mK\mM\mB~+\hspace{2cm}
\]
\[\hspace{4cm}+~\mK\mG\mB\mM-\mG\mK\mB\mM+\mG\mM\mK\mB-\mG\mK\mM\mB-\mG\mB\mM\mK+\mG\mB\mK\mM~\grh~=~0\,.
\]
\  \eindebewijs

\ \\[1mm]

\newpage
\renewcommand\thesection{\Alph{section}}\setcounter{section}{0} 

\section{ Addendum on Free Gauge Fields}
If we put 
\[
\sG_A(\osx)=\sG(~\ldots, \aF_{\mu\nu}(\osx),\ldots;\ldots, \aF^\dagger_{\th\rho}(\osx),\ldots;\osx\,)~=~\hspace{6cm}
\]
\[
 \  \hspace{1cm} =\sG(\ldots,\pa_\mu\cA_\nu-\pa_\nu\cA_\mu-\llh\cA_\mu\lk\cA_\nu\rlh,\ldots\,;
\ldots,\pa_\mu\cA_\nu^\da-\pa_\nu\cA_\mu^\da+\llh\cA_\mu^\da\lk\cA_\nu^\da\rlh,\ldots\,;\osx\,)~=~~~~~~~~~~
\]
\be ~\hspace{2cm}=~\sL(\os{\aA}(\osx)\,;\, \os{\aA}^\dagger(\osx)\,;\,\ldots,\pa_\mu\os{\aA}(\osx),\ldots\,;\,\ldots,\pa_\mu\os{\aA}^\dagger(\osx),\ldots\,;\,\osx)\,,
\ee
with $ \os{\aA} = \col[\ldots,\aA_\mu,\ldots]$, which now plays the role of $\dpsi$ in section 2, we get, in accordance with our notation in section 2,
\be
\begin{array}{ccccccccccl}
\sL^{\dg{(o)}}_A &=& \text{row }[&\ldots&\ldots&\ldots&\ldots&-\sum_{\mu=1}^N\llh\hat{\sG}^{(\dg{\mu\ka})}_A \lk \cA_\mu\rlh\,&\ldots&\ldots &\hspace{-2mm}] \\[2mm]
\sL^{\dg{(1)}}_A &=& \text{row }[& 0 &\sG^{(\dg{12})}_A  & \sG^{(\dg{13})}_A & \ldots & \sG^{(\dg{1\ka})}_A  &  \ldots  &  \sG^{(\dg{1N})}_A &\hspace{-2mm}]\\[2mm]
\sL^{\dg{(2)}}_A &=& \text{row }[&-\sG^{(\dg{12})}_A  &0 &\sG^{(\dg{23})}_A & \ldots & \sG^{(\dg{2\ka})}_A  &  \ldots  &  \sG^{(\dg{2N})}_A&\hspace{-2mm}]\\[2mm]
\sL^{\dg{(3)}}_A &=& \text{row }[&-\sG^{(\dg{13})}_A  &-\sG^{(\dg{23})}_A &   0   & \ldots & \sG^{(\dg{3\ka})}_A  &  \ldots  &  \sG^{(\dg{3N})}_A  &\hspace{-2mm}]\\[2mm]
\ldots &=& \text{row }[&\ldots&\ldots &   \ldots  & \ldots & \ldots  &  \ldots  &  \ldots  &\hspace{-2mm}]\\[2mm]
\sL^{\dg{(\ka)}}_A &=& \text{row }[&-\sG^{(\dg{1\ka})}_A  &-\sG^{(\dg{2\ka})}_A & -\sG^{(\dg{3\ka})}_A   & \ldots &    0     &  \ldots  &  \sG^{(\dg{\ka N})}_A  &\hspace{-2mm}]\\[2mm]
\ldots &=& \text{row }[&\ldots&\ldots &   \ldots  & \ldots & \ldots  &  \ldots  &  \ldots  &\hspace{-2mm}]\\[2mm]
\sL^{\dg{(N)}}_A &=& \text{row }[&-\sG^{(\dg{1N})}_A  &-\sG^{(\dg{2N})}_A & -\sG^{(\dg{3N})}_A   & \ldots &-\sG^{(\dg{\ka N})}_A     &  \ldots  &  0  &\hspace{-2mm}]\\[2mm]
\end{array}
\ee
With convention (3.15) the lower $N$ rows of this table simplify to
\be
\sL^{\dg{(\mu)}}_A ~=~ \text{row }[\ldots,\hat{\sG}^{(\dg{\mu\ka})}_A\,,\ldots]\,,~~~1\leq\mu,\ka\leq N\,.
\ee
Table (A.2) enables to reduce the proof of Theorem 3.2 to an application of Theorem 2.4.\\[1mm]
Because of property (3.7) it is obvious that all 'components' of $\sL^{\dg{(\mu\star)}}_A\,,0\leq\mu\leq N,$ are the {\em hermitean transposed} of the components of $\sL^{\dg{(\mu)}}_A\,,0\leq\mu\leq N\,$. Only for $\sL^{\dg{(o\star)}}_A$ this is not immediately obvious. Let us check it  in an {\em ad hoc} way by
calculating the $\ka$-th component of $\sL^{\dg{(o\star)}}_A$. In (A.1) replace $\llh\cA_\mu^\da\lk\cA_\nu^\da\rlh$ by the perturbation $\llh\cA_\mu^\da+\vep\del_{\mu\ka} H\lk\cA_\nu^\da+\vep\del_{\nu\ka} H\rlh$. Now differentiate the result to $\vep$. At $\vep=0$ it becomes
\[
\sum_{1\leq\mu<\nu\leq N}\Tr\glh\sG^{(\dg{\mu\nu\star})}_A:\llh\del_{\mu\ka} H\lk\aA_\nu^\da\rlh+\llh\aA_\mu^\da\lk\del_{\nu\ka} H\rlh\grh~= \hspace{6cm}
\]
\[
 \ \hspace{2cm}  =~\sum_{\ka<\nu\leq N}\Tr\glh\sG^{(\dg{\ka\nu\star})}_A:\llh H\lk\aA_\nu^\da\rlh\grh~+~
\sum_{1\leq\mu<\ka}\Tr\glh\sG^{(\dg{\mu\ka\star})}_A:\llh \aA_\mu^\da\lk H\rlh\grh~=
\]
\[
=~\sum_{\ka<\nu\leq N}\Tr\glh\llh\aA_\nu^\da\lk\sG^{(\dg{\ka\nu\star})}_A\rlh: H\grh~+~
\sum_{1\leq\mu<\ka}\Tr\glh\llh\sG^{(\dg{\mu\ka\star})}_A\lk\aA_\mu^\da\rlh: H\grh~=~
\Tr\glh\sum_{\mu=1}^N\llh\hat{\sG}^{(\dg{\mu\ka\star})}_A\lk\aA_\mu^\da\rlh: H\grh\,.
\]
Finally one finds
\[
\glh\sum_{\mu=1}^N\llh\hat{\sG}^{(\dg{\mu\ka\star})}_A\lk\aA_\mu^\da\rlh\grh^\da=-\sum_{\mu=1}^N\llh\hat{\sG}^{(\dg{\mu\ka})}_A \lk \cA_\mu\rlh\,.
\]
\ \\[0.5mm]
{\bf Remark on Thm 4.9-b:} If it happens that
\[
\sG(\ldots,e^{s S_\mu^\la}\pa_\la\cA_\nu-e^{s S_\nu^\th}\pa_\th\cA_\mu-\llh\cA_\mu\lk\cA_\nu\rlh,\ldots\,;
\ldots,e^{s S_\mu^\la}\pa_\la\cA_\nu^\da-e^{s S_\nu^\th}\pa_\th\cA_\mu^\da+\llh\cA_\mu^\da\lk\cA_\nu^\da\rlh,\ldots\,;\osx\,)~=~
\]
\[
~\hspace{1cm} =\sG(\ldots,\pa_\mu\cA_\nu-\pa_\nu\cA_\mu-\llh\cA_\mu\lk\cA_\nu\rlh,\ldots\,;
\ldots,\pa_\mu\cA_\nu^\da-\pa_\nu\cA_\mu^\da+\llh\cA_\mu^\da\lk\cA_\nu^\da\rlh,\ldots\,;\osx\,)~+~\ms{O}(s^2)\,,
\]
it follows that
\[
\Re\sum_{\mu<\nu}\Tr\glh\sG^{(\dg{\mu\nu})}_A:S_\mu^\la\pa_\la\aA_\nu-S_\nu^\th\pa_\th\aA_\mu\grh~=~0\,.
\]

\section{ Electromagnetism}

Some more details on Example 3.4B:

\[
\sG_A=\sum_{0\leq\mu<\nu\leq3}(-1)^{\del_{\mu0}+\del_{\nu0}}\Tr\big[\aF_{\mu\nu}^\dagger  \aF_{\mu\nu} \,\big]\,
\]
\[
\sG^{(\dg{01})}_A=-\aF_{01}^\dagger~~~\sG^{(\dg{02})}_A=-\aF_{02}^\dagger~~~\sG^{(\dg{03})}_A=-\aF_{03}^\dagger~~~\sG^{(\dg{12})}_A=\aF_{12}^\dagger~~~\sG^{(\dg{13})}_A=\aF_{13}^\dagger
~~~\sG^{(\dg{23})}_A=\aF_{23}^\dagger
\]
Now (3.19) reads, for $0\leq\ka\leq3$,
\[
\begin{array}{cl}
\ka=0\,:  &  \pa_1\sG^{(\dg{01})}_A+\pa_2\sG^{(\dg{02})}_A+\pa_3\sG^{(\dg{03})}_A=\\[2mm]
          &~~~~ =-\pa_1(\pa_0\aA^\dagger_1-\pa_1\aA^\dagger_0)-\pa_2(\pa_0\aA^\dagger_2-\pa_2\aA^\dagger_0)-\pa_3(\pa_0\aA^\dagger_3-\pa_3\aA^\dagger_0)\\[2mm]
          &~~~~ =-\pa_0(\pa_1\aA^\dagger_1+\pa_2\aA^\dagger_2+\pa_3\aA^\dagger_3)+\pa_1\pa_1\aA^\dagger_0+\pa_2\pa_2\aA^\dagger_0+\pa_3\pa_3\aA^\dagger_0\\[2mm]
\ka=1\,:  &  -\pa_0\sG^{(\dg{01})}_A+\pa_2\sG^{(\dg{12})}_A+\pa_3\sG^{(\dg{13})}_A=\\
          &~~~~ =\pa_0(\pa_0\aA^\dagger_1-\pa_1\aA^\dagger_0)+\pa_2(\pa_1\aA^\dagger_2-\pa_2\aA^\dagger_1)+\pa_3(\pa_1\aA^\dagger_3-\pa_3\aA^\dagger_1)\\[2mm]
          &~~~~ =\pa_0\pa_0\aA^\dagger_1+\pa_1(-\pa_0\aA^\dagger_0+\pa_1\aA^\dagger_1+\pa_2\aA^\dagger_2+\pa_3\aA^\dagger_3)-(\pa_1\pa_1+\pa_2\pa_2+\pa_3\pa_3)\aA^\dagger_1\\[2mm]
\ka=2\,:  &  -\pa_0\sG^{(\dg{02})}_A-\pa_1\sG^{(\dg{12})}_A+\pa_3\sG^{(\dg{23})}_A=\\
          &~~~~ =\pa_0(\pa_0\aA^\dagger_2-\pa_2\aA^\dagger_0)-\pa_1(\pa_1\aA^\dagger_2-\pa_2\aA^\dagger_1)+\pa_3(\pa_2\aA^\dagger_3-\pa_3\aA^\dagger_2)\\[2mm]
          &~~~~ =\pa_0\pa_0\aA^\dagger_2+\pa_2(-\pa_0\aA^\dagger_0+\pa_1\aA^\dagger_1+\pa_2\aA^\dagger_2+\pa_3\aA^\dagger_3)-(\pa_1\pa_1+\pa_2\pa_2+\pa_3\pa_3)\aA^\dagger_2\\[2mm]
\ka=3\,:  &  -\pa_0\sG^{(\dg{03})}_A-\pa_1\sG^{(\dg{13})}_A-\pa_2\sG^{(\dg{23})}_A=\\
          &~~~~ =\pa_0(\pa_0\aA^\dagger_3-\pa_3\aA^\dagger_0)-\pa_1(\pa_1\aA^\dagger_3-\pa_3\aA^\dagger_1)-\pa_2(\pa_2\aA^\dagger_3-\pa_3\aA^\dagger_2)\\[2mm]
          &~~~~ =\pa_0\pa_0\aA^\dagger_3+\pa_3(-\pa_0\aA^\dagger_0+\pa_1\aA^\dagger_1+\pa_2\aA^\dagger_2+\pa_3\aA^\dagger_3)-(\pa_1\pa_1+\pa_2\pa_2+\pa_3\pa_3)\aA^\dagger_3\\[2mm]                             
\end{array}
\]
If we put $\aA_0^\dagger=-\Phi$ and $\col[\aA_1^\dagger\,,\aA_2^\dagger\,,\aA_3^\dagger]=\os{A}$ we get Maxwell's equations 'in potential form'
\be
\left\{\begin{array}{rcc}
\paf{}{t}\div\os{A} +\Delta \Phi  & =& 0\\[2mm]
\paft{}{t}\os{A}-\Delta\os{A}+\grad\big(\paf{}{t}\Phi+\div\os{A}\big) &=&\os{0}
\end{array}
\right.
\ee
If the pair $\os{A},\os{B}$ satisfies this pair, then the pair $\os{E}=-\paf{\os{A}}{t}-\grad\Phi\,,\,\os{B}=\rot\os{A}\,$, satisfies the classical homogeneous Maxwell equations:
\[
\pa_t\os{B}=\rot\,\pa_t\os{A}=\rot(-\os{E}-\grad\Phi)=-\rot\os{E}
\]
\[
\pa_t\os{E}=\pa_t\pa_t\os{A}-\grad\,\pa_t\Phi=-\Delta\os{A}+\grad\,\div\os{A}=\rot\,\rot\os{A}=\rot\os{B}
\]

Finally, imposing the 'Lorenz-Gauge' $~\paf{}{t}\Phi+\div\os{A}=0$, we find the usual wave equations for $\Phi$ and \os{A}.

Any solution to the system (B.1) can be reduced to a solution which satisfies the Lorentz condition, by means of a 'gauge transform'
$\Phi\mapsto \Phi-\pa_t\Lambda\,,\,\os{A}\mapsto\os{A}-\grad\Lambda$, leading to the same $\os{E},\os{B}$-fields. cf. Jackson [J], p.241.

\ \\[2mm]Similar results can be found for more general free fields governed by

\[\sG_1=\sum_{\mu\nu\ths\rhos}\,g^{\mu\ths}g^{\nu\rhos}\Tr\glh J\cF^\dagger_{\ths\rhos}J^{-1}\cF_{\mu\nu}\grh\,.
\]
\section*{References}
\begin{itemize}
\item[{\bf [DM]}] W. Drechsler, M.E. Mayor: Fiber Bundle Techniques in Gauge Theories. Lecture Notes in Physics 67. Springer Berlin 1977.
\item[{\bf [H]}] L. H{\"o}rmander: Complex Analysis in Several Variables. North Holland Publ. Co. 1973.
\item[{\bf [J]}] J.D. Jackson: Classical Electrodynamics 3rd ed. John Wiley. N.J. 1998.
\item[{\bf [JP]}] E.M. de Jager, H.G.J. Pijls: Proc. Sem. Mathematical Structures in Field Theories 1981-1982. ISBN 90 6196 2781. CWI Amsterdam 1984.
\item[{\bf [M]}] A. Messiah: Quantum Mechanics Vol II. North Holland Publ. Co. 1962.
\item[{\bf [W]}] S. Weinberg: The Quantum Theory of Fields. Vol I: CUP, Cambridge 1995.
\end{itemize}

\end{document}